\documentclass[10pt,12pt]{article}%
\usepackage{amsmath}
\usepackage{amsfonts}
\usepackage{amssymb}

\providecommand{\U}[1]{\protect\rule{.1in}{.1in}}

\begin{document}

\title{A spectral inequality for degenerated operators \\and applications}
\author{R\'{e}mi Buffe\thanks{Universit\'e de Lorraine, CNRS, Inria, IECL, F-54000 Nancy, France. E-mail address: remi.buffe@inria.fr} , Kim Dang Phung\thanks{
Institut Denis Poisson, CNRS, UMR 7013,Universit\'{e} d'Orl\'{e}ans, BP
6759,45067 Orl\'{e}ans Cedex 2, France. E-mail address:
kim\_dang\_phung@yahoo.fr. The author was partially supported by the
R\'{e}gion Centre (France) - THESPEGE Project.}}
\date{15/03/2018 \ }
\maketitle

\bigskip

Abstract .- In this paper we establish a Lebeau-Robbiano spectral inequality
for a degenerated one dimensional elliptic operator and show how it can be
used to impulse control and finite time stabilization for a degenerated
parabolic equation.

\bigskip

R\'{e}sum\'{e} .- Dans cet article, on s'int\'{e}r\`{e}sse \`{a}
l'in\'{e}galit\'{e} spectrale de type Lebeau-Robbiano sur la somme de
fonctions propres pour une famille d'op\'{e}rateurs d\'{e}g\'{e}n\'{e}r\'{e}s.
Les applications sont donn\'{e}es en th\'{e}orie du contr\^{o}le comme le
contr\^{o}le impulsionnel et la stabilisation en temps fini.

\bigskip

\bigskip

\section{Introduction and main results}

\bigskip

The purpose of this article is to prove spectral properties for a family of
degenerate operators acting on the interval $\left(  0,1\right)  $. We shall
consider linear operators ${\mathcal{P}}$ in $L^{2}\left(  0,1\right)  $,
defined by
\[
\left\{
\begin{array}
[c]{ll}%
{\mathcal{P}}=-\frac{d}{dx}\left(  x^{\alpha}\frac{d}{dx}\right)  \text{ ,
with }\alpha\in(0,2)\text{ ,} & \\
D({\mathcal{P}})=\left\{  \vartheta\in H_{\alpha}^{1}\left(  0,1\right)
\text{; }{\mathcal{P}}\vartheta\in L^{2}(0,1)\text{ and BC}_{\alpha}%
(\vartheta)=0\right\}  \text{ ,} &
\end{array}
\right.
\]
where%
\[
H_{\alpha}^{1}(0,1):=\left\{  \vartheta\in L^{2}(0,1)\text{;}\;\vartheta\text{
is absolutely continuous in }(0,1)\text{,}\;\int_{0}^{1}x^{\alpha}%
|\vartheta^{\prime}|^{2}<\infty\text{, }\vartheta(1)=0\right\}  \text{ ,}%
\]
and
\[
\text{BC}_{\alpha}(\vartheta)=%
\begin{cases}
\vartheta_{|_{x=0}}\text{ ,} & \text{for }\alpha\in\lbrack0,1)\text{ ,}\\
(x^{\alpha}\vartheta^{\prime})_{|_{x=0}}\text{ ,} & \text{for }\alpha
\in\lbrack1,2)\text{ .}%
\end{cases}
\]
We remind that ${\mathcal{P}}$ is a closed self-adjoint positive densely
defined operator, with compact resolvent. As a consequence, the following
spectral decomposition holds: There exists a countable family of
eigenfunctions $\Phi_{j}$ associated with eigenvalues $\lambda_{j}$ such that

\begin{itemize}
\item $\left\{  \Phi_{j}\right\}  _{j\geq1}$ forms an Hilbert basis of
$L^{2}(0,1)$

\item ${\mathcal{P}}\Phi_{j}=\lambda_{j}\Phi_{j}$

\item $0<\lambda_{1}\leq\lambda_{2}\leq\dots\leq\lambda_{k}\rightarrow+\infty$ .
\end{itemize}

An explicit expression of the eigenvalues is given in \cite{Gu} for the the
weakly degenerate case $\alpha\in\left(  0,1\right)  $, and in \cite{Mo} for
the strongly degenerate case $\alpha\in\lbrack1,2)$, and depends on the Bessel
function of first kind (see \cite{MM}). Also, we have the following asymptotic
formula: $\lambda_{k}\sim C\left(  \alpha\right)  k^{2}$ as $k\rightarrow
\infty$.

\bigskip

We are interested on the spectral inequality for the sum of eigenfunctions.
Our main result is as follows.

\bigskip

Theorem 1.1 .- \textit{Let }$\omega$\textit{ be an open and nonempty subset of
}$(0,1)$\textit{. There exist constants }$C>0$\textit{ and }$\sigma\in\left(
0,1\right)  $\textit{ such that }%
\[
\sum_{\lambda_{j}\leq\Lambda}|a_{j}|^{2}\leq Ce^{C\Lambda^{\sigma}}%
\int_{\omega}\left\vert \sum_{\lambda_{j}\leq\Lambda}a_{j}\Phi_{j}\right\vert
^{2}\text{ ,}%
\]
\textit{for all }$\left\{  a_{j}\right\}  \in\mathbb{R}$\textit{ and }%
$\Lambda>0$\textit{. Further, }%
\[
\sigma=%
\begin{cases}
3/4\text{ ,} & \text{\textit{if} }\alpha\in\left(  0,2\right)  \left\backslash
\left\{  1\right\}  \right.  \text{ ,}\\
3/\left(  2\gamma\right)  \text{ \textit{for any} }\gamma\in\left(
0,2\right)  \text{ ,} & \text{\textit{if} }\alpha=1\text{ .}%
\end{cases}
\]

\bigskip

Two different kinds of approach have been developed to obtain the spectral
inequality for the sum of eigenfunctions: A first one is due to Lebeau and
Robbiano \cite{LR} and is based on a Carleman estimate for an elliptic
operator, whereas a second one appears in a remark in \cite{AEWZ} and is based
on an observation estimate at one point in time for a parabolic equation. Note
that in the standard setting of uniformly elliptic operator, $\sigma=1/2$ (see
\cite{L}, \cite{JL}, \cite{LZ}, \cite{Lu}, \cite{Mi}, \cite{LRR1},
\cite{LRLR}). In the present paper we will establish a new Carleman estimate
for an associated degenerated elliptic operator. Because of the degeneracy of
the coefficients of the operator ${\mathcal{P}}$, we make use of a new weight
function in the design of the Carleman estimate. The subtle difference between
the cases $\alpha\in\left(  0,2\right)  \left\backslash \left\{  1\right\}
\right.  $ and $\alpha=1$ is related to the existence of a Hardy type
inequality for the $H_{\alpha}^{1}$ norm. Indeed, for $\alpha=1$, the desired
Hardy inequality fails to hold.

\bigskip

Many applications to such spectral inequality have been developed, in
particular in control theory (see \cite{L}, \cite{LZ}, \cite{BN}, \cite{Le},
\cite{LRM}, \cite{BPS}). Let $\omega$ be an open and nonempty subset of
$(0,1)$ and denote $1_{\omega}$ the characteristic function of a given
subdomain $\omega$. We present the following two results.

\bigskip

Theorem 1.2 .- \textit{Let }$E\subset\left(  0,T\right)  $\textit{ be a
measurable set of positive measure. For all }$y^{0}\in L^{2}(0,1)$\textit{,
there exists }$f\in L^{2}\left(  \omega\times E\right)  $\textit{ such that
the solution }$y=y\left(  x,t\right)  $\textit{ of }%
\[
\left\{
\begin{array}
[c]{ll}%
\partial_{t}y-\partial_{x}\left(  x^{\alpha}\partial_{x}y\right)
=1_{\omega\times E}f\text{ ,} & \text{\textit{in} }(0,1)\times(0,T)\text{ ,}\\
\text{BC}_{\alpha}(y)=0\text{ ,} & \text{\textit{on} }(0,T)\text{ ,}\\
y_{|_{x=1}}=0\text{ ,} & \text{\textit{on} }(0,T)\text{ ,}\\
y_{|_{t=0}}=y^{0}\text{ ,} & \text{\textit{in} }(0,1)\text{ ,}%
\end{array}
\right.
\]
\textit{satisfies }$y(\cdot,T)=0$\textit{.}

\bigskip

Theorem 1.3 .- \textit{There is }$\left(  t_{m}\right)  _{m\in\mathbb{N}}%
$\textit{ a increasing sequence of positive real numbers converging to }%
$T>0$\textit{ and }$\left(  \mathcal{F}_{m}\right)  _{m\in\mathbb{N}}$\textit{
a sequence of linear bounded operators from }$L^{2}(0,1)$\textit{ into }%
$L^{2}(0,1)$\textit{ such that for any }$z_{0}\in L^{2}(0,1)$\textit{, the
solution }$z=z\left(  x,t\right)  $\textit{ to }%
\[
\left\{
\begin{array}
[c]{ll}%
\partial_{t}z-\partial_{x}\left(  x^{\alpha}\partial_{x}z\right)
=\sum\limits_{m\in\mathbb{N}}\delta_{t=\left(  t_{m+1}+t_{m}\right)
/2}\otimes1_{\omega}\mathcal{F}_{m}\left(  z_{|_{t=t_{m}}}\right)  \text{ ,} &
\text{\textit{in} }(0,1)\times(0,T)\text{ ,}\\
\text{BC}_{\alpha}(z)=0\text{ ,} & \text{\textit{on} }(0,T)\text{ ,}\\
z_{|_{x=1}}=0\text{ ,} & \text{\textit{on} }(0,T)\text{ ,}\\
z_{|_{t=0}}=z_{0}\text{ ,} & \text{\textit{in} }(0,1)\text{ ,}%
\end{array}
\right.
\]
\textit{satisfies }$\underset{t\rightarrow T_{-}}{\text{lim}}\left\Vert
z\left(  \cdot,t\right)  \right\Vert _{L^{2}\left(  \Omega\right)  }%
=0$\textit{.}

\bigskip

Here $\delta_{t=\left(  t_{m+1}+t_{m}\right)  /2}$ denotes the Dirac measure
at $t=\left(  t_{m+1}+t_{m}\right)  /2$. Note that the above system
equivalently reads%
\[
\left\{
\begin{array}
[c]{ll}%
\partial_{t}z-\partial_{x}\left(  x^{\alpha}\partial_{x}z\right)  =0\text{ ,}
& \text{\textit{for} }t\in\mathbb{R}^{+}\backslash\bigcup\limits_{m\geq
0}\left(  \frac{t_{m}+t_{m+1}}{2}\right)  \ \text{,}\\
z\left(  \cdot,\frac{t_{m}+t_{m+1}}{2}\right)  =z\left(  \cdot,\left(
\frac{t_{m}+t_{m+1}}{2}\right)  _{-}\right)  +1_{\omega}\mathcal{F}_{m}\left(
z\left(  \cdot,t_{m}\right)  \right)  \text{ ,} & \text{\textit{for any
integer} }m\geq0\text{ ,}\\
\text{BC}_{\alpha}(z)=0\text{ ,} & \text{\textit{on} }(0,T)\text{ ,}\\
z_{|_{x=1}}=0\text{ ,} & \text{\textit{on} }(0,T)\text{ ,}\\
z_{|_{t=0}}=z_{0}\text{ ,} & \text{\textit{in} }(0,1)\text{ .}%
\end{array}
\right.
\]
Theorem 3.1 is new approach to steer the solution to zero at time $T$ and can
be seen as a finite time stabilization for the degenerated heat equation by
impulse control. This can be compared with \cite{CN}. The standard
null-controllability problem is given when $E=\left(  0,T\right)  $ and has
been studied in \cite{CMV}. It is now well-known that the null controllability
for higher degeneracies ($\alpha\geq2$) fails to hold (see \cite{CMV2} and the
references therein). We also refer to \cite{ABCF}, where the
null-controllability result has been extended to more general degeneracies at
the boundary. When the control is located at the boundary where the degeneracy
occurs, we refer to \cite{Gu,CTY,MRR}. We finally refer to the recent book
\cite{CMV2} and the references therein for a full description of the field.
Note that an estimation of the cost of controllability for small $T>0$, as
well as for $\alpha\rightarrow2^{-}$ has been recently obtained in \cite{CMV3}.

\bigskip

The outline of the paper is as follows. In Section 2, we present the key
inequalities needed to prove Theorem 1.1 as Hardy inequality and Carleman
inequality. Section 3 is devoted to obtaining the applications of the spectral
inequality in control theory as observation estimates, impulse approximate
controllability, null controllability on measurable set in time (see Theorem
3.4) and finite time stabilization (see Theorem 3.5). Theorem 1.2 and Theorem
1.3 are direct consequence of Theorem 3.4 and Theorem 3.5 respectively.

\bigskip

\bigskip

\section{Key inequalities}

\bigskip

This section is devoted to the statement of the key inequalities: Hardy
inequality and Carleman inequality, that will enable us to prove Theorem 1.1.
The proof of the Carleman inequality is given at the end of this section.

\bigskip

\subsection{Hardy inequality and boundary conditions}

\bigskip

The following Hardy inequality shall play a central role in what follows. The
proof can be found in \cite{CMV}, \cite{OK}.

\bigskip

Lemma 2.1 .- \textit{Let }$\vartheta$\textit{ be a locally absolutely
continuous functions on }$(0,1)$\textit{ such that }$\displaystyle\int_{0}%
^{1}x^{\alpha}|\vartheta^{\prime}|^{2}<\infty$\textit{. Then we have }%
\[
\int_{0}^{1}x^{\alpha-2}|\vartheta|^{2}\leq\frac{4}{(2-\alpha)^{2}}\int
_{0}^{1}x^{\alpha}|\vartheta^{\prime}|^{2}\text{ ,}%
\]
\textit{if one of the following assumption holds:}%
\[%
\begin{array}
[c]{cc}%
i) & \alpha\in\left(  0,1\right)  \text{ \textit{and} }\vartheta_{|_{x=0}%
}=0\text{ ,}\\
ii) & \alpha\in(1,2)\text{ \textit{and} }\vartheta_{|_{x=1}}=0\text{ .}%
\end{array}
\]

\bigskip

We also have the following lemma, that shall be useful when estimating the
boundary terms arising from integration by parts in the strongly degenerate
case $\alpha\in\lbrack1,2)$. The proof can be found in \cite{CMV}.

\bigskip

Lemma 2.2 .- \textit{Let }$\alpha\in\lbrack1,2)$\textit{ and }$\vartheta\in
H_{\alpha}^{1}(0,1)$\textit{. Then }$(x|\vartheta|^{2})_{|_{x=0}}=0$\textit{.}

\bigskip

\bigskip

\subsection{Global Carleman estimate near the degeneracy}

\bigskip

In this section, we shall state the crucial tool, i.e. a global Carleman
estimate near the degeneracy of an elliptic operator.

\bigskip

Introduce, for $S_{0}>s_{0}>0$,
\[
Z=(-S_{0},S_{0})\times(0,1)\text{ ,}\quad Y=(-s_{0},s_{0})\times(0,1)\text{ .}%
\]
First, we shall write
\begin{equation}
{\mathcal{Q}}:=-\partial_{s}^{2}+{\mathcal{P}}=-\partial_{s}^{2}-\partial
_{x}\left(  x^{\alpha}\partial_{x}\right)  \text{ ,} \tag{2.2.1}\label{2.2.1}%
\end{equation}
here $\left(  s,x\right)  \in Z$. The weight function we choose is of the
form
\begin{equation}
\varphi(\tau,s,x)=\tau\frac{x^{2-\alpha}}{2-\alpha}-\frac{\tau^{\gamma/3}}%
{\nu}s^{2}\text{ ,} \tag{2.2.2}\label{2.2.2}%
\end{equation}
where $\tau,\nu>0$ are two large parameters,
\begin{equation}%
\begin{cases}
\gamma=2\text{ ,} & \text{for }\alpha\in\left(  0,2\right)  \left\backslash
\left\{  1\right\}  \right.  \text{ ,}\\
\gamma<2\text{ ,} & \text{for }\alpha=1\text{ ,}%
\end{cases}
\tag{2.2.3}\label{2.2.3}%
\end{equation}
and with $\nu$ fixed sufficiently large. Note that this weight function is
completely decoupled in the two directions, in particular with respect to the
dependency in $\tau$. In the case $\alpha=1$, the Hardy inequality in Lemma
2.1 does not hold, and this is the reason of our subtle choice of weight
(\ref{2.2.2}). Next, we shall set
\[
{\mathcal{Q}}_{\varphi}:=e^{\varphi}{\mathcal{Q}}e^{-\varphi}\text{ .}%
\]
Finally, we state a global estimate for functions of $C^{\infty}((-S_{0}%
,S_{0}),D({\mathcal{P}}))$, with the proper weight function $\varphi$ given by
(\ref{2.2.2}) to handle the degeneracy at $x=0$.

\bigskip

Theorem 2.1 .- \textit{There exist }$\tau_{0}>0$\textit{, and }$\nu_{0}%
>0$\textit{ such that for }$\gamma>0$\textit{ defined in (\ref{2.2.3}), there
exists }$c>0$\textit{ such that }%
\[
\tau^{\gamma}||v||_{L^{2}(Z)}^{2}+\tau\int_{Z}x^{\alpha}|\partial_{x}%
v|^{2}+\tau^{3}\int_{Z}x^{2-\alpha}|v|^{2}+\mathcal{B}(v)\leq c||{\mathcal{Q}%
}_{\varphi}v||_{L^{2}(Z)}^{2}\text{ ,}%
\]
\textit{for all }$\tau\geq\tau_{0}$\textit{, and for all }$v\in C^{\infty
}((-S_{0},S_{0}),D(P))$\textit{, where }$\mathcal{B}$\textit{ is a quadratic
form satisfying}%
\[%
\begin{array}
[c]{ll}%
\frac{1}{2}\mathcal{B}(v)\geq & -\tau\displaystyle\int_{-S_{0}}^{S_{0}%
}|\partial_{x}v_{|_{x=1}}|^{2}+2\displaystyle\frac{\tau^{\gamma/3}}{\nu}%
\int_{0}^{1}\left[  s|\partial_{s}v|^{2}\right]  _{s=-S_{0}}^{s=S_{0}%
}+2\displaystyle\frac{\tau^{\gamma/3}}{\nu}\int_{0}^{1}\left[  v\partial
_{s}v\right]  _{s=-S_{0}}^{s=S_{0}}\\
& +8\displaystyle\frac{\tau^{\gamma}}{\nu^{3}}\int_{0}^{1}\left[  s^{3}%
|v|^{2}\right]  _{s=-S_{0}}^{s=S_{0}}-2\tau\int_{0}^{1}\left[  x\partial
_{s}v\partial_{x}v\right]  _{s=-S_{0}}^{s=S_{0}}-\tau\int_{0}^{1}\left[
v\partial_{s}v\right]  _{s=-S_{0}}^{s=S_{0}}\\
& -2\displaystyle\frac{\tau^{\gamma/3}}{\nu_{0}}\int_{0}^{1}\left[  x^{\alpha
}s|\partial_{x}v|^{2}\right]  _{s=-S_{0}}^{s=S_{0}}+2\displaystyle\frac
{\tau^{2+\gamma/3}}{\nu_{0}}\int_{0}^{1}\left[  x^{2-\alpha}s|v|^{2}\right]
_{s=-S_{0}}^{s=S_{0}}\text{ .}%
\end{array}
\]

\bigskip

Note that in the above Theorem 2.1, boundary conditions are prescribed through
the membership in the domain of $\mathcal{P}$. The proof will be given at the
end of this section.

\bigskip

In \cite{CMV}, the authors established a parabolic Carleman estimate for a
class of degenerated operators, in the spirit of \cite{FI}, with a weight
linked to geodesic distance to the singularity $\{x=0\}$, that is a weight of
the form
\begin{equation}
\widetilde{\varphi}(x,t)=\frac{x^{2-\alpha}-1}{(t(T-t))^{4}}\text{ .}
\tag{2.2.4}\label{2.2.4}%
\end{equation}
In the present article, the design of the weight function $\varphi$ is similar
to (\ref{2.2.4}). However, as we have to deal with an additional variable $s$
(see the operator (\ref{2.2.1})), we also weaken the weight function in the
$s$ direction (see the weight (\ref{2.2.2}) which is anisotropic with respect
to powers of the Carleman large parameter $\tau$).

\bigskip

\bigskip

\subsection{Inequality with weight for a specific sum of eigenfunctions}

\bigskip

A classical trick on quantitative uniqueness consists on transferring
properties for elliptic equation into an estimate for parabolic operator (see
\cite{Li}). Here, we naturally reproduce this idea for the sum of
eigenfunctions (see \cite{L}, \cite{JL}, \cite{LR}, \cite{LZ}, \cite{CSL},
\cite{Lu}, \cite{LRL}, \cite{Le}).

\bigskip

We define the following function space, depending on the frequency parameter
$\Lambda\geq1$,
\[
{\mathcal{X}}_{\Lambda}:=\left\{  \mathbf{u}(s,x)=\sum_{\lambda_{j}\leq
\Lambda}\frac{\text{sinh}(\sqrt{\lambda_{j}}(s+S_{0}))}{\sqrt{\lambda_{j}}%
}a_{j}\Phi_{j}(x)\text{;}\;a_{j}\in{\mathbb{R}}\right\}  \text{ .}%
\]
We then go back to a weighted estimate for functions $\mathbf{u}%
\in{\mathcal{X}}_{\Lambda}$. Notice that ${\mathcal{Q}}\mathbf{u}=0$.

\bigskip

Corollary 2.1 .- \textit{Let }$\gamma>0$\textit{ defined in (\ref{2.2.3}).
There exist }$\tau_{0}>0$\textit{ and }$c>0$\textit{ such that }%
\[
\tau^{\gamma}\left\Vert e^{\varphi}\mathbf{u}\right\Vert _{L^{2}(Z)}^{2}%
+\tau\int_{Z}x^{\alpha}|e^{\varphi}\partial_{x}\mathbf{u}|^{2}+\tau^{3}%
\int_{Z}x^{2-\alpha}|e^{\varphi}\mathbf{u}|^{2}\leq c\tau\int_{-S_{0}}^{S_{0}%
}|e^{\varphi}\partial_{x}\mathbf{u}_{|_{x=1}}|^{2}\text{ ,}%
\]
\textit{for }$\tau=\tau_{0}\Lambda^{3/\left(  2\gamma\right)  }$\textit{, and
for all }$\mathbf{u}\in{\mathcal{X}}_{\Lambda}$\textit{, }$\Lambda\geq
1$\textit{.}

\bigskip

Proof .- We shall apply the Carleman estimate in Theorem 2.1 to $v=e^{\varphi
}\mathbf{u}$, with $\mathbf{u}\in{\mathcal{X}}_{\Lambda}$. Recall that
${\mathcal{Q}}\mathbf{u}=0$. Clearly, we have $v_{|_{s=-S_{0}}}=\partial
_{x}v_{|_{s=-S_{0}}}=0$, and also ${\mathcal{Q}}_{\varphi}v=0$. By Theorem
2.1, this yields
\begin{equation}
\tau^{\gamma}\left\Vert v\right\Vert _{L^{2}(Z)}^{2}+\tau\int_{Z}x^{\alpha
}|\partial_{x}v|^{2}+\tau^{3}\int_{Z}x^{2-\alpha}|v|^{2}+\mathcal{B}%
(v)\leq0\text{ ,} \tag{2.3.1}\label{2.3.1}%
\end{equation}
with%
\begin{equation}%
\begin{array}
[c]{ll}%
\frac{1}{2}\mathcal{B}(v)\geq & -\tau\displaystyle\int_{-S_{0}}^{S_{0}%
}|\partial_{x}v_{|_{x=1}}|^{2}+2\displaystyle\frac{\tau^{\gamma/3}}{\nu}%
\int_{0}^{1}S_{0}|\partial_{s}v_{|_{s=S_{0}}}|^{2}+2\displaystyle\frac
{\tau^{\gamma/3}}{\nu}\int_{0}^{1}v_{|_{s=S_{0}}}\partial_{s}v_{|_{s=S_{0}}}\\
& +8\displaystyle\frac{\tau^{\gamma}}{\nu^{3}}\int_{0}^{1}S_{0}^{3}%
|v_{|_{s=S_{0}}}|^{2}-2\tau\displaystyle\int_{0}^{1}x\partial_{s}%
v_{|_{s=S_{0}}}\partial_{x}v_{|_{s=S_{0}}}-\tau\displaystyle\int_{0}%
^{1}v_{|_{s=S_{0}}}\partial_{s}v_{|_{s=S_{0}}}\\
& -2\displaystyle\frac{\tau^{\gamma/3}}{\nu_{0}}\int_{0}^{1}x^{\alpha}%
S_{0}|\partial_{x}v_{|_{s=S_{0}}}|^{2}+2\displaystyle\frac{\tau^{2+\gamma/3}%
}{\nu_{0}}\int_{0}^{1}x^{2-\alpha}S_{0}|v_{|_{s=S_{0}}}|^{2}\text{ .}%
\end{array}
\tag{2.3.2}\label{2.3.2}%
\end{equation}
We first work with volumic terms (from now, the notation $A\lesssim B$ means
that there exists a constant $c>0$, independent on the concerned parameters
such that $A\leq cB$). We have
\begin{equation}
\tau\int_{Z}x^{\alpha}|e^{\varphi}\partial_{x}\mathbf{u}|^{2}\lesssim\tau
\int_{Z}x^{\alpha}|\partial_{x}v|^{2}+\tau\int_{Z}x^{\alpha}|(\partial
_{x}\varphi)v|^{2}\lesssim\tau\int_{Z}x^{\alpha}|\partial_{x}v|^{2}+\tau
^{3}\int_{Z}x^{2-\alpha}|v|^{2}\text{ .} \tag{2.3.3}\label{2.3.3}%
\end{equation}
Therefore, from (\ref{2.3.1}), there exists $c>0$ such that
\begin{equation}
\tau^{\gamma}\left\Vert e^{\varphi}\mathbf{u}\right\Vert _{L^{2}(Z)}^{2}%
+\tau\int_{Z}x^{\alpha}|e^{\varphi}\partial_{x}\mathbf{u}|^{2}+\tau^{3}%
\int_{Z}x^{2-\alpha}|e^{\varphi}\mathbf{u}|^{2}+c\mathcal{B}(v)\leq0\text{ .}
\tag{2.3.4}\label{2.3.4}%
\end{equation}
Remark that
\begin{equation}
\tau^{\gamma}\left\Vert v\right\Vert _{L^{2}(Z)}^{2}\geq\tau^{\gamma
}\left\Vert v\right\Vert _{L^{2}(Y)}^{2}\gtrsim\tau^{\gamma}e^{-\frac{2}{\nu
}\tau^{\gamma/3}s_{0}^{2}}\int_{0}^{1}|e^{\tau\frac{x^{2-\alpha}}{2-\alpha}%
}\sum\limits_{\lambda_{j}\leq\Lambda}a_{j}\Phi_{j}|^{2}\text{ ,}
\tag{2.3.5}\label{2.3.5}%
\end{equation}
and then using (\ref{2.3.5}) with (\ref{2.3.4}), we obtain that there exists
$c>0$ such that
\begin{equation}%
\begin{array}
[c]{ll}
& \quad e^{-\frac{2}{\nu}\tau^{\gamma/3}s_{0}^{2}}\displaystyle\int_{0}%
^{1}|e^{\tau\frac{x^{2-\alpha}}{2-\alpha}}\sum\limits_{\lambda_{j}\leq\Lambda
}a_{j}\Phi_{j}(x)|^{2}\\
& +\tau^{\gamma}||e^{\varphi}\mathbf{u}||_{L^{2}(Z)}^{2}+\tau\int_{Z}%
x^{\alpha}|e^{\varphi}\partial_{x}\mathbf{u}|^{2}+\tau^{3}\int_{Z}x^{2-\alpha
}|e^{\varphi}\mathbf{u}|^{2}+c\mathcal{B}(v)\leq0\text{ .}%
\end{array}
\tag{2.3.6}\label{2.3.6}%
\end{equation}
Second, we handle the boundary terms. We have, using Young inequality,
\[
\tau\left\vert \int_{0}^{1}x\partial_{s}v_{|_{s=S_{0}}}\partial_{x}%
v_{|_{s=S_{0}}}\right\vert \lesssim\tau^{\gamma/3}\int_{0}^{1}x^{\alpha
}|\partial_{x}v_{|_{s=S_{0}}}|^{2}+\tau^{2\gamma/3}\int_{0}^{1}x^{2-\alpha
}|\partial_{s}v_{|_{s=S_{0}}}|^{2}\text{ .}%
\]
Note that from the form of $v$, we obtain that%
\[%
\begin{array}
[c]{ll}
& \quad\tau^{\gamma/3}\displaystyle\int_{0}^{1}x^{\alpha}|\partial
_{x}v_{|_{s=S_{0}}}|^{2}\\
& \lesssim\tau^{\gamma}\displaystyle\int_{0}^{1}x^{2-\alpha}|v_{|_{s=S_{0}}%
}|^{2}+\tau^{\gamma/3}\displaystyle\int_{0}^{1}x^{\alpha}|e^{\varphi}%
\partial_{x}\mathbf{u}_{|_{s=S_{0}}}|^{2}\\
& \lesssim\tau^{\gamma}\displaystyle\int_{0}^{1}x^{2-\alpha}|v_{|_{s=S_{0}}%
}|^{2}+\tau^{\gamma/3}\Lambda\displaystyle\int_{0}^{1}|v_{|_{s=S_{0}}}|^{2}\\
& \lesssim e^{-\frac{2}{\nu}\tau^{\gamma/3}S_{0}^{2}}e^{cS_{0}\sqrt{\Lambda}%
}\left(  \tau^{\gamma}\displaystyle\int_{0}^{1}x^{2-\alpha}|e^{\tau
\frac{x^{2-\alpha}}{2-\alpha}}\sum\limits_{\lambda_{j}\leq\Lambda}a_{j}%
\Phi_{j}|^{2}+\tau^{\gamma/3}\Lambda\displaystyle\int_{0}^{1}|e^{\tau
\frac{x^{2-\alpha}}{2-\alpha}}\sum\limits_{\lambda_{j}\leq\Lambda}a_{j}%
\Phi_{j}|^{2}\right)  \text{ .}%
\end{array}
\]
Taking $\tau=\tau_{0}\Lambda^{3/\left(  2\gamma\right)  }$, with $\tau_{0}>0$
sufficiently large, yields
\[
\tau^{\gamma/3}\int_{0}^{1}x^{\alpha}|\partial_{x}v_{|_{s=S_{0}}}|^{2}%
\lesssim\tau_{0}^{\gamma}\Lambda^{3/2}e^{-\frac{2}{\nu}\tau_{0}^{\gamma
/3}\sqrt{\Lambda}S_{0}^{2}}e^{cS_{0}\sqrt{\Lambda}}\int_{0}^{1}|e^{\tau
\frac{x^{2-\alpha}}{2-\alpha}}\sum\limits_{\lambda_{j}\leq\Lambda}a_{j}%
\Phi_{j}|^{2}\text{ .}%
\]
Also, using the form of $v$, and then taking $\tau=\tau_{0}\Lambda^{3/\left(
2\gamma\right)  }$, with $\tau_{0}>0$ sufficiently large, one can deduce that%
\[%
\begin{array}
[c]{ll}%
\tau^{2\gamma/3}\displaystyle\int_{0}^{1}x^{2-\alpha}|\partial_{s}%
v_{|_{s=S_{0}}}|^{2} & \lesssim\tau^{\gamma}\displaystyle\int_{0}%
^{1}x^{2-\alpha}|v_{|_{s=S_{0}}}|^{2}+\tau^{2\gamma/3}\displaystyle\int
_{0}^{1}x^{2-\alpha}|e^{\varphi}\partial_{s}\mathbf{u}_{|_{s=S_{0}}}|^{2}\\
& \lesssim e^{-\frac{2}{\nu}\tau^{\gamma/3}S_{0}^{2}}e^{cS_{0}\sqrt{\Lambda}%
}\left(  \tau^{\gamma}+\tau^{2\gamma/3}\Lambda\right)  \displaystyle\int
_{0}^{1}x^{2-\alpha}|e^{\tau\frac{x^{2-\alpha}}{2-\alpha}}\sum\limits_{\lambda
_{j}\leq\Lambda}a_{j}\Phi_{j}|^{2}\\
& \lesssim\tau_{0}^{\gamma}\Lambda^{2}e^{-\frac{2}{\nu}\tau_{0}^{\gamma
/3}\sqrt{\Lambda}S_{0}^{2}}e^{cS_{0}\sqrt{\Lambda}}\displaystyle\int_{0}%
^{1}|e^{\tau\frac{x^{2-\alpha}}{2-\alpha}}\sum\limits_{\lambda_{j}\leq\Lambda
}a_{j}\Phi_{j}|^{2}\text{ .}%
\end{array}
\]
Using the same arguments, taking $\tau=\tau_{0}\Lambda^{3/\left(
2\gamma\right)  }$, with $\tau_{0}>0$ sufficiently large, yields
\[
\tau\left\vert \int_{0}^{1}v_{|_{s=S_{0}}}\partial_{s}v_{|_{s=S_{0}}%
}\right\vert \lesssim\tau_{0}^{\gamma}\Lambda^{2}e^{-\frac{2}{\nu}\tau
_{0}^{\gamma/3}\sqrt{\Lambda}S_{0}^{2}}e^{cS_{0}\sqrt{\Lambda}}\int_{0}%
^{1}|e^{\tau\frac{x^{2-\alpha}}{2-\alpha}}\sum\limits_{\lambda_{j}\leq\Lambda
}a_{j}\Phi_{j}|^{2}\text{ .}%
\]
At this point, we see that all the negative terms at $s=S_{0}$ in
(\ref{2.3.2}) are bounded by
\begin{equation}
\tau_{0}^{\gamma}\Lambda^{2}e^{-\frac{2}{\nu}\tau_{0}^{\gamma/3}\sqrt{\Lambda
}S_{0}^{2}}e^{cS_{0}\sqrt{\Lambda}}\int_{0}^{1}|e^{\tau\frac{x^{2-\alpha}%
}{2-\alpha}}\sum\limits_{\lambda_{j}\leq\Lambda}a_{j}\Phi_{j}|^{2}\text{ .}
\tag{2.3.7}\label{2.3.7}%
\end{equation}
As a result, since $s_{0}<S_{0}$, the quantity (\ref{2.3.7}) can be dominated
by the first term in the left hand side of (\ref{2.3.6}), by taking $\tau
=\tau_{0}\Lambda^{3/\left(  2\gamma\right)  }$, with $\tau_{0}$ sufficiently
large. Hence, from the boundary terms $\mathcal{B}(v)$, it only remains
$-\tau\int_{-S_{0}}^{S_{0}}|\partial_{x}v_{|_{x=1}}|^{2}$. But,
\[
-\tau\int_{-S_{0}}^{S_{0}}|\partial_{x}v_{|_{x=1}}|^{2}=-\tau\int_{-S_{0}%
}^{S_{0}}|e^{\varphi}\partial_{x}\mathbf{u}_{|_{x=1}}|^{2}\text{ ,}%
\]
by using the boundary conditions. This ends the proof. \hfill$\Box$

\bigskip

\bigskip

\subsection{Proof of Theorems 1.1 and 2.1}

\bigskip

\subsubsection{Proof of the spectral inequality}

\bigskip

This section is devoted to proving Theorem 1.2. First we establish the
spectral inequality with an observation at the boundary $\{x=1\}$ by applying
Corollary 2.1, we have%
\[
\tau^{\gamma}\left\Vert e^{\varphi}\mathbf{u}\right\Vert _{L^{2}(Z)}^{2}%
+\tau\int_{Z}x^{\alpha}|e^{\varphi}\partial_{x}\mathbf{u}|^{2}+\tau^{3}%
\int_{Z}x^{2-\alpha}|e^{\varphi}\mathbf{u}|^{2}\leq c\tau\int_{-S_{0}}^{S_{0}%
}|e^{\varphi}\partial_{x}\mathbf{u}_{|_{x=1}}|^{2}\text{ ,}%
\]
with for all $\tau=\tau_{0}\Lambda^{3/\left(  2\gamma\right)  }$, and for all
$u\in X_{\Lambda}$, $\Lambda\geq1$. Arguing as in (\ref{2.3.5}), we obtain%
\[
\tau^{\gamma}\left\Vert e^{\varphi}\mathbf{u}\right\Vert _{L^{2}(Z)}^{2}%
\geq\tau^{\gamma}\left\Vert e^{\varphi}\mathbf{u}\right\Vert _{L^{2}(Y)}%
^{2}\gtrsim\tau^{\gamma}e^{-\frac{2}{\nu}\tau^{\gamma/3}s_{0}^{2}}\int_{0}%
^{1}|e^{\tau\frac{x^{2-\alpha}}{2-\alpha}}\sum\limits_{\lambda_{j}\leq\Lambda
}a_{j}\Phi_{j}|^{2}\text{ .}%
\]
Bounding the weight functions, and keeping in mind that $\tau=\tau_{0}%
\Lambda^{3/\left(  2\gamma\right)  }$, one can deduce that there exist
$C_{1},C_{2},C_{3}>0$ such that
\[
e^{-C_{1}\sqrt{\Lambda}}\left\Vert \sum\limits_{\lambda_{j}\leq\Lambda}%
a_{j}\Phi_{j}\right\Vert _{L^{2}(0,1)}^{2}\leq C_{2}e^{C_{3}\Lambda
^{3/(2\gamma)}}\left\vert \sum\limits_{\lambda_{j}\leq\Lambda}a_{j}\Phi
_{j}^{\prime}(1)\right\vert ^{2}\text{ ,}%
\]
which is the spectral inequality with a boundary observation. Then, as the
region $\{x=1\}$ is away from the singularity, the operator ${\mathcal{P}}$ is
uniformly elliptic there, and therefore it is classical (see for instance
\cite{R}, \cite{LRL}, \cite{L}) that we can propagate the observation
$\{x=1\}$ to $\{s=-S_{0}\}\times\omega$ by using classical Carleman estimates
to obtain the desired spectral inequality.\hfill$\Box$

\bigskip

\subsubsection{Proof of Theorem 2.1}

\bigskip

Here, we give the proof of the global Carleman estimate near the degeneracy in
Theorem 2.1.

\bigskip

Recall that ${\mathcal{Q}}_{\varphi}=e^{\varphi}{\mathcal{Q}}e^{-\varphi}$ and
therefore
\[%
\begin{array}
[c]{ll}%
{\mathcal{Q}}_{\varphi} & =-(\partial_{s}-(\partial_{s}\varphi))^{2}%
+(\partial_{x}-(\partial_{x}\varphi))x^{\alpha}(\partial_{x}-(\partial
_{x}\varphi))\\
& =-\partial_{s}^{2}-|\partial_{s}\varphi|^{2}+2(\partial_{s}\varphi
)\partial_{s}+\partial_{s}^{2}\varphi+{\mathcal{P}}-x^{2\alpha}|\partial
_{x}\varphi|^{2}+2x^{\alpha}(\partial_{x}\varphi)\partial_{x}+\partial
_{x}\left(  x^{\alpha}\partial_{x}\varphi\right)  \text{ .}%
\end{array}
\]
Now, we decompose ${\mathcal{Q}}_{\varphi}$ into four parts:
\[
{\mathcal{Q}}_{\varphi}={\mathcal{S}}_{x}+{\mathcal{S}}_{s}+{\mathcal{A}}%
_{x}+{\mathcal{A}}_{s}\text{ ,}%
\]
where ${\mathcal{S}}_{x}+{\mathcal{S}}_{s}$ is the symmetric part and
${\mathcal{A}}_{x}+{\mathcal{A}}_{s}$ is the skew-symmetric part of the full
conjugated operator. Using the definition of the weight function
(\ref{2.2.2}), we have
\[
{\mathcal{S}}_{x}={\mathcal{P}}-\tau^{2}x^{2-\alpha}\text{ ,}\quad
{\mathcal{S}}_{s}=-\partial_{s}^{2}-4\frac{\tau^{2\gamma/3}}{\nu^{2}}%
s^{2}\text{ ,}\quad{\mathcal{A}}_{x}=2\tau x\partial_{x}+\tau\text{ ,}%
\quad{\mathcal{A}}_{s}=-4\frac{\tau^{\gamma/3}}{\nu}s\partial_{s}-2\frac
{\tau^{\gamma/3}}{\nu}\text{ .}%
\]
Let $v\in C^{\infty}((-S_{0},S_{0}),D({\mathcal{P}}))$. We begin by noting
that%
\[%
\begin{array}
[c]{ll}%
\left\Vert {\mathcal{Q}}_{\varphi}v\right\Vert _{L^{2}(Z)}^{2} &
=||{\mathcal{S}}v||_{L^{2}(Z)}^{2}+||{\mathcal{A}}v||_{L^{2}(Z)}^{2}+2\left(
{\mathcal{S}}v,{\mathcal{A}}v\right)  _{Z}\\
& \geq||{\mathcal{S}}v||_{L^{2}(Z)}^{2}+2\left[  \left(  {\mathcal{S}}%
_{x}v,{\mathcal{A}}_{x}v\right)  _{Z}+\left(  {\mathcal{S}}_{s}v,{\mathcal{A}%
}_{s}v\right)  _{Z}+\left(  {\mathcal{S}}_{x}v,{\mathcal{A}}_{s}v\right)
_{Z}+\left(  {\mathcal{S}}_{s}v,{\mathcal{A}}_{x}v\right)  _{Z}\right]  \text{
.}%
\end{array}
\]
The proof is divided into three steps. Each step corresponds to the
computation of one of the above scalar products.

\bigskip

\textit{First Step. We begin with the first scalar product }$\left(
{\mathcal{S}}_{x}v,{\mathcal{A}}_{x}v\right)  _{Z}$\textit{.}

Lemma 2.3 .- \textit{We have }%
\begin{equation}
\left(  {\mathcal{S}}_{x}v,{\mathcal{A}}_{x}v\right)  _{Z}=\tau(2-\alpha
)\int_{Z}x^{\alpha}|\partial_{x}v|^{2}+\tau^{3}(2-\alpha)\int_{Z}x^{2-\alpha
}|v|^{2}+\mathcal{B}_{0}(v)\text{ ,} \tag{4.1}\label{4.1}%
\end{equation}
\textit{with }%
\[
\mathcal{B}_{0}(v)=-\tau\int_{-S_{0}}^{S_{0}}\left[  x^{\alpha+1}|\partial
_{x}v|^{2}\right]  _{x=0}^{x=1}-\tau\int_{-S_{0}}^{S_{0}}\left[  x^{\alpha
}v\partial_{x}v\right]  _{x=0}^{x=1}-\tau^{3}\int_{-S_{0}}^{S_{0}}\left[
x^{3-\alpha}|v|^{2}\right]  _{x=0}^{x=1}\text{ .}%
\]

\bigskip

The proof of this lemma will be provided later. Using the Hardy inequality of
Lemma 2.1 in (\ref{4.1}), there exists $c>0$ such that for all $\alpha
\in\left(  0,2\right)  \left\backslash \left\{  1\right\}  \right.  $,
\[
c\left(  {\mathcal{S}}_{x}v,{\mathcal{A}}_{x}v\right)  _{Z}\geq\tau^{2}%
\int_{Z}|v|^{2}+\tau(2-\alpha)\int_{Z}x^{\alpha}|\partial_{x}v|^{2}+\tau
^{3}(2-\alpha)\int_{Z}x^{2-\alpha}|v|^{2}+\mathcal{B}_{0}(v)\text{ .}%
\]
In the particular case $\alpha=1$, using the Hardy inequality, for all
$\alpha^{\prime}\in(1,2)$, there exists $c^{\prime}>0$ such that
\begin{equation}
\tau\int_{Z}x|\partial_{x}v|^{2}\geq\tau\int_{Z}x^{\alpha^{\prime}}%
|\partial_{x}v|^{2}\geq c^{\prime}\tau\int_{Z}x^{\alpha^{\prime}-2}%
|v|^{2}\text{ .} \tag{4.2}\label{4.2}%
\end{equation}
As a result, interpolating (\ref{4.2}) with (\ref{4.1}), for all $\gamma
\in(0,2)$, there exists $c^{\prime}>0$ such that
\begin{equation}
c^{\prime}\left(  {\mathcal{S}}_{x}v,{\mathcal{A}}_{x}v\right)  _{Z}\geq
\tau^{\gamma}\int_{Z}|v|^{2}+\tau(2-\alpha)\int_{Z}x^{\alpha}|\partial
_{x}v|^{2}+\tau^{3}(2-\alpha)\int_{Z}x^{2-\alpha}|v|^{2}+\mathcal{B}%
_{0}(v)\text{ .} \tag{4.3}\label{4.3}%
\end{equation}
Hence, (\ref{4.3}) holds for all $\alpha\in\left(  0,2\right)  $, with
$\gamma$ defined in (\ref{2.2.3}). We now focus on boundary terms
$\mathcal{B}_{0}$. We have, using boundary conditions described in $H_{\alpha
}^{1}(0,1)$ and Lemma 2.2,
\[
\mathcal{B}_{0}(v)=-\tau\int_{-S_{0}}^{S_{0}}|\partial_{x}v_{|_{x=1}}%
|^{2}\text{ .}%
\]

\textit{Second step. We then compute the second scalar product }$\left(
{\mathcal{S}}_{s}v,{\mathcal{A}}_{s}v\right)  _{Z}$\textit{.}

Lemma 2.4 .- \textit{We have }%
\[
\left(  {\mathcal{S}}_{s}v,{\mathcal{A}}_{s}v\right)  _{Z}=-4\frac
{\tau^{\gamma/3}}{\nu}\int_{Z}|\partial_{s}v|^{2}-16\frac{\tau^{\gamma}}%
{\nu^{3}}\int_{Z}s^{2}|v|^{2}+\mathcal{B}_{1}(v)\text{ ,}%
\]
\textit{with }%
\[
\mathcal{B}_{1}(v):=2\frac{\tau^{\gamma/3}}{\nu}\int_{0}^{1}\left[
s|\partial_{s}v|^{2}\right]  _{s=-S_{0}}^{s=S_{0}}+2\frac{\tau^{\gamma/3}}%
{\nu}\int_{0}^{1}\left[  v\partial_{s}v\right]  _{s=-S_{0}}^{s=S_{0}}%
+8\frac{\tau^{\gamma}}{\nu^{3}}\int_{0}^{1}\left[  s^{3}|v|^{2}\right]
_{s=-S_{0}}^{s=S_{0}}\text{ .}%
\]

\bigskip

The proof of this lemma will be provided later. The two volumic terms in Lemma
2.4 are non-positive, and need a particular attention. We set
\[
K_{1}(v):=-4\frac{\tau^{\gamma/3}}{\nu}\int_{Z}|\partial_{s}v|^{2}\text{
,}\quad K_{2}(v):=-16\frac{\tau^{\gamma}}{\nu^{3}}\int_{Z}s^{2}|v|^{2}\text{
.}%
\]
Introduce ${\mathcal{S}}={\mathcal{S}}_{x}+{\mathcal{S}}_{s}$, we have the
following relation
\[
-\int_{Z}|\partial_{s}v|^{2}=\int_{Z}({\mathcal{S}}v)v-4\frac{\tau^{2\gamma
/3}}{\nu^{2}}\int_{Z}s^{2}|v|^{2}+\int_{Z}({\mathcal{P}}v)v-\tau^{2}\int
_{Z}x^{2-\alpha}|v|^{2}\text{ ,}%
\]
and we then deduce%
\[%
\begin{array}
[c]{ll}%
K_{1}(v) & =-4\displaystyle\frac{\tau^{\gamma/3}}{\nu}\int_{Z}({\mathcal{S}%
}v)v-16\displaystyle\frac{\tau^{\gamma}}{\nu^{3}}\int_{Z}s^{2}|v|^{2}%
+4\displaystyle\frac{\tau^{\gamma/3}}{\nu}\int_{Z}({\mathcal{P}}%
v)v-4\displaystyle\frac{\tau^{2+\gamma/3}}{\nu}\int_{Z}x^{2-\alpha}|v|^{2}\\
& =-4\displaystyle\frac{\tau^{\gamma/3}}{\nu}\int_{Z}({\mathcal{S}}%
v)v+K_{2}(v)+4\displaystyle\frac{\tau^{\gamma/3}}{\nu}\int_{Z}({\mathcal{P}%
}v)v-4\displaystyle\frac{\tau^{2+\gamma/3}}{\nu}\int_{Z}x^{2-\alpha}%
|v|^{2}\text{ .}%
\end{array}
\]
As a result, using integration by parts and Young inequality,%
\[%
\begin{array}
[c]{ll}
& \quad\left(  {\mathcal{S}}_{s}v,{\mathcal{A}}_{s}v\right)  _{Z}\\
& =-4\displaystyle\frac{\tau^{\gamma/3}}{\nu}\int_{Z}({\mathcal{S}}%
v)v+K_{2}(v)+4\displaystyle\frac{\tau^{\gamma/3}}{\nu}\int_{Z}({\mathcal{P}%
}v)v-4\displaystyle\frac{\tau^{2+\gamma/3}}{\nu}\int_{Z}x^{2-\alpha}%
|v|^{2}+\mathcal{B}_{1}(v)\\
& \geq-\displaystyle\frac{2}{\nu}||{\mathcal{S}}v||_{L^{2}(Z)}^{2}%
-2\displaystyle\frac{\tau^{2\gamma/3}}{\nu}||v||_{L^{2}(Z)}^{2}+K_{2}%
(v)+4\displaystyle\frac{\tau^{\gamma/3}}{\nu}\int_{Z}x^{\alpha}|\partial
_{x}v|^{2}-4\displaystyle\frac{\tau^{2+\gamma/3}}{\nu}\int_{Z}x^{2-\alpha
}|v|^{2}\\
& \quad+\widetilde{\mathcal{B}}_{1}(v)\text{ ,}%
\end{array}
\]
with
\[
\widetilde{\mathcal{B}}_{1}(v)=\mathcal{B}_{1}(v)-\frac{4\tau^{\gamma/3}}{\nu
}\int_{-S_{0}}^{S_{0}}\left[  x^{\alpha}v\partial_{x}v\right]  _{x=0}%
^{x=1}\text{ .}%
\]
Note that, using boundary conditions, we have $\widetilde{\mathcal{B}}%
_{1}=\mathcal{B}_{1}$. Summing up, fixing $\nu:=\nu_{0}>0$ sufficiently large,
and taking $\tau\geq\tau_{0}$, with $\tau_{0}>0$ sufficiently large, there
exists $c>0$ such that%
\[%
\begin{array}
[c]{ll}
& \quad c\left(  \left\Vert {\mathcal{S}}v\right\Vert _{L^{2}(Z)}^{2}+2\left(
{\mathcal{S}}_{x}v,{\mathcal{A}}_{x}v\right)  +2\left(  {\mathcal{S}}%
_{s}v,{\mathcal{A}}_{s}v\right)  \right) \\
& \geq\tau^{\gamma}\left\Vert v\right\Vert _{L^{2}(Z)}^{2}+\tau
\displaystyle\int_{Z}x^{\alpha}|\partial_{x}v|^{2}+\tau^{3}\displaystyle\int
_{Z}x^{2-\alpha}|v|^{2}+2\mathcal{B}_{0}(v)+2\mathcal{B}_{1}(v)\text{ .}%
\end{array}
\]

\textit{Third step. It remains to estimate the crossed-terms }$\left(
{\mathcal{S}}_{x}v,{\mathcal{A}}_{s}v\right)  _{Z}+\left(  {\mathcal{S}}%
_{s}v,{\mathcal{A}}_{x}v\right)  _{Z}$\textit{.}

Lemma 2.5 .- \textit{We have, on the one hand}%
\[%
\begin{array}
[c]{ll}%
\left(  {\mathcal{S}}_{s}v,{\mathcal{A}}_{x}v\right)  _{Z}=\mathcal{B}%
_{2}(v):= & \tau\displaystyle\int_{-S_{0}}^{S_{0}}\left[  x|\partial_{s}%
v|^{2}\right]  _{x=0}^{x=1}-2\tau\displaystyle\int_{0}^{1}\left[
x\partial_{s}v\partial_{x}v\right]  _{s=-S_{0}}^{s=S_{0}}\\
& -\tau\displaystyle\int_{0}^{1}\left[  v\partial_{s}v\right]  _{s=-S_{0}%
}^{s=S_{0}}-4\displaystyle\frac{\tau^{1+2\gamma/3}}{\nu_{0}^{2}}\int_{-S_{0}%
}^{S_{0}}\left[  s^{2}x|v|^{2}\right]  _{x=0}^{x=1}\text{ ,}%
\end{array}
\]
\textit{and on the other hand}%
\[%
\begin{array}
[c]{ll}%
\left(  {\mathcal{S}}_{x}v,{\mathcal{A}}_{s}v\right)  _{Z}=\mathcal{B}%
_{3}(v):= & 4\displaystyle\frac{\tau^{\gamma/3}}{\nu_{0}}\int_{-S_{0}}^{S_{0}%
}\left[  x^{\alpha}s\partial_{x}v\partial_{s}v\right]  _{x=0}^{x=1}%
-2\displaystyle\frac{\tau^{\gamma/3}}{\nu_{0}}\int_{0}^{1}\left[  x^{\alpha
}s|\partial_{x}v|^{2}\right]  _{s=-S_{0}}^{s=S_{0}}\\
& +2\displaystyle\frac{\tau^{\gamma/3}}{\nu_{0}}\int_{-S_{0}}^{S_{0}}\left[
x^{\alpha}v\partial_{x}v\right]  _{x=0}^{x=1}+2\displaystyle\frac
{\tau^{2+\gamma/3}}{\nu_{0}}\int_{0}^{1}\left[  x^{2-\alpha}s|v|^{2}\right]
_{s=-S_{0}}^{s=S_{0}}\text{ .}%
\end{array}
\]

\bigskip

The proof of this lemma will be provided later. Note that using boundary
conditions given in $H_{\alpha}^{1}(0,1)$ as well as Lemma 2.2, we have
\[
\mathcal{B}_{2}(v)\geq-2\tau\int_{0}^{1}\left[  x\partial_{s}v\partial
_{x}v\right]  _{s=-S_{0}}^{s=S_{0}}-\tau\int_{0}^{1}\left[  v\partial
_{s}v\right]  _{s=-S_{0}}^{s=S_{0}}\text{ ,}%
\]
and
\[
\mathcal{B}_{3}(v)=-2\frac{\tau^{\gamma/3}}{\nu_{0}}\int_{0}^{1}\left[
x^{\alpha}s|\partial_{x}v|^{2}\right]  _{s=-S_{0}}^{s=S_{0}}+2\frac
{\tau^{2+\gamma/3}}{\nu_{0}}\int_{0}^{1}\left[  x^{2-\alpha}s|v|^{2}\right]
_{s=-S_{0}}^{s=S_{0}}\text{ .}%
\]
Now setting $\mathcal{B}=2\left(  \mathcal{B}_{0}+\mathcal{B}_{1}%
+\mathcal{B}_{2}+\mathcal{B}_{3}\right)  $ yields the sought result.
\hfill$\Box$

\bigskip

\subsubsection{Proof of Lemma 2.3}

We recall that
\[
{\mathcal{S}}_{x}={\mathcal{P}}-\tau^{2}x^{2-\alpha}\text{ ,}\quad
{\mathcal{A}}_{x}=2\tau x\partial_{x}+\tau\text{ .}%
\]
We shall denote by ${\mathcal{I}}_{ij}$ the scalar product between the
$i^{th}$ term of ${\mathcal{S}}_{x}$ with the $j^{th}$ term of ${\mathcal{A}%
}_{x}$. Let us compute first%
\[%
\begin{array}
[c]{ll}%
{\mathcal{I}}_{11} & =-2\tau\displaystyle\int_{Z}\partial_{x}x^{\alpha
}\partial_{x}vx\partial_{x}v\\
& =2\tau\displaystyle\int_{Z}x^{\alpha}|\partial_{x}v|^{2}+\tau
\displaystyle\int_{Z}x^{1+\alpha}\partial_{x}\left(  |\partial_{x}%
v|^{2}\right)  -2\tau\displaystyle\int_{-S_{0}}^{S_{0}}\left[  x^{1+\alpha
}|\partial_{x}v|^{2}\right]  _{x=0}^{x=1}\\
& =(1-\alpha)\tau\displaystyle\int_{Z}x^{\alpha}|\partial_{x}v|^{2}%
-\tau\displaystyle\int_{-S_{0}}^{S_{0}}\left[  x^{1+\alpha}|\partial_{x}%
v|^{2}\right]  _{x=0}^{x=1}\text{ .}%
\end{array}
\]
Second, we have
\[
{\mathcal{I}}_{12}=-\tau\int_{Z}v\partial_{x}x^{\alpha}\partial_{x}v=\tau
\int_{Z}x^{\alpha}|\partial_{x}v|^{2}-\tau\int_{-S_{0}}^{S_{0}}\left[
x^{\alpha}v\partial_{x}v\right]  _{x=0}^{x=1}\text{ .}%
\]
Third, we see that
\[
{\mathcal{I}}_{21}=-2\tau^{3}\int_{Z}x^{3-\alpha}v\partial_{x}v=-\tau^{3}%
\int_{Z}x^{3-\alpha}\partial_{x}\left(  |v|^{2}\right)  =(3-\alpha)\tau
^{3}\int_{Z}x^{2-\alpha}|v|^{2}-\tau^{3}\int_{-S_{0}}^{S_{0}}\left[
x^{3-\alpha}|v|^{2}\right]  _{x=0}^{x=1}\text{ .}%
\]
Finally, we can check that
\[
{\mathcal{I}}_{22}=-\tau^{3}\int_{Z}x^{2-\alpha}|v|^{2}\text{ ,}%
\]
and we end the proof of Lemma 2.3 by summing the above four quantities.\hfill
$\Box$

\bigskip

\subsubsection{Proof of Lemma 2.4}

We recall that
\[
{\mathcal{S}}_{s}=-\partial_{s}^{2}-4\frac{\tau^{2\gamma/3}}{\nu^{2}}%
s^{2}\text{ ,}\quad{\mathcal{A}}_{s}=-4\frac{\tau^{\gamma/3}}{\nu}%
s\partial_{s}-2\frac{\tau^{\gamma/3}}{\nu}\text{ .}%
\]
We shall denote by ${\mathcal{I}}_{ij}$ the scalar product between the
$i^{th}$ term of ${\mathcal{S}}_{s}$ with the $j^{th}$ term of ${\mathcal{A}%
}_{s}$. Let us compute the ${\mathcal{I}}_{ij}$, $1\leq i,j\leq2$, by
integrations by parts
\[
{\mathcal{I}}_{11}=\frac{4\tau^{\gamma/3}}{\nu}\int_{Z}s\partial_{s}%
^{2}vs\partial_{s}v=\frac{2\tau^{\gamma/3}}{\nu}\int_{Z}s\partial_{s}\left(
|\partial_{s}v|^{2}\right)  =-\frac{2\tau^{\gamma/3}}{\nu}\int_{Z}%
|\partial_{s}v|^{2}+\frac{2\tau^{\gamma/3}}{\nu}\int_{0}^{1}\left[
s|\partial_{s}v|^{2}\right]  _{s=-S_{0}}^{s=S_{0}}\text{ ,}%
\]%
\[
{\mathcal{I}}_{12}=\frac{2\tau^{\gamma/3}}{\nu}\int_{Z}v\partial_{s}%
^{2}v=-\frac{2\tau^{\gamma/3}}{\nu}\int_{Z}|\partial_{s}v|^{2}+\frac
{2\tau^{\gamma/3}}{\nu}\int_{0}^{1}\left[  v\partial_{s}v\right]  _{s=-S_{0}%
}^{s=S_{0}}\text{ ,}%
\]%
\[
{\mathcal{I}}_{21}=\frac{16\tau^{\gamma}}{\nu^{3}}\int_{Z}s^{3}v\partial
_{s}v=8\frac{\tau^{\gamma}}{\nu^{3}}\int_{Z}s^{3}\partial_{s}\left(
|v|^{2}\right)  =-\frac{24\tau^{\gamma}}{\nu^{3}}\int_{Z}s^{2}|v|^{2}%
+\frac{8\tau^{\gamma}}{\nu^{3}}\int_{0}^{1}\left[  s^{3}|v|^{2}\right]
_{s=-S_{0}}^{s=S_{0}}\text{ ,}%
\]
and
\[
{\mathcal{I}}_{22}=8\frac{\tau^{\gamma}}{\nu^{3}}\int_{Z}s^{2}|v|^{2}\text{ .}%
\]
Summing all the ${\mathcal{I}}_{ij}$ yields the sought result of Lemma
2.4.\hfill$\Box$

\bigskip

\subsubsection{Proof of Lemma 2.5}

We recall that
\[
{\mathcal{S}}_{x}={\mathcal{P}}-\tau^{2}x^{2-\alpha}\text{ ,}\quad
{\mathcal{S}}_{s}=-\partial_{s}^{2}-4\frac{\tau^{2\gamma/3}}{\nu^{2}}%
s^{2}\text{ ,}\quad{\mathcal{A}}_{x}=2\tau x\partial_{x}+\tau\text{ ,}%
\quad{\mathcal{A}}_{s}=-4\frac{\tau^{\gamma/3}}{\nu}s\partial_{s}-2\frac
{\tau^{\gamma/3}}{\nu}\text{ .}%
\]
We first compute the scalar product $\left(  {\mathcal{S}}_{s}v,{\mathcal{A}%
}_{x}v\right)  _{Z}$. We shall denote by ${\mathcal{I}}_{ij}$ the scalar
product between the $i^{th}$ term of ${\mathcal{S}}_{s}$ with the $j^{th}$
term of ${\mathcal{A}}_{x}$. We have%
\[%
\begin{array}
[c]{ll}%
{\mathcal{I}}_{11} & =-2\tau\displaystyle\int_{Z}x\partial_{s}^{2}%
v\partial_{x}v=\tau\displaystyle\int_{Z}x\partial_{x}\left(  |\partial
_{s}v|^{2}\right)  -2\tau\displaystyle\int_{0}^{1}\left[  x\partial
_{s}v\partial_{x}v\right]  _{s=-S_{0}}^{s=S_{0}}\\
& =-\tau\displaystyle\int_{Z}|\partial_{s}v|^{2}+\tau\displaystyle\int
_{-S_{0}}^{S_{0}}\left[  x|\partial_{s}v|^{2}\right]  _{x=0}^{x=1}%
-2\tau\displaystyle\int_{0}^{1}\left[  x\partial_{s}v\partial_{x}v\right]
_{s=-S_{0}}^{s=S_{0}}\text{ ,}%
\end{array}
\]%
\[
{\mathcal{I}}_{12}=-\tau\int_{Z}v\partial_{s}^{2}v=\tau\int_{Z}|\partial
_{s}v|^{2}-\tau\int_{0}^{1}\left[  v\partial_{s}v\right]  _{s=-S_{0}}%
^{s=S_{0}}\text{ ,}%
\]%
\[%
\begin{array}
[c]{ll}%
{\mathcal{I}}_{21} & =-\displaystyle\frac{8\tau^{2\gamma/3+1}}{\nu^{2}}%
\int_{Z}s^{2}xv\partial_{x}v=-\displaystyle\frac{4\tau^{2\gamma/3+1}}{\nu^{2}%
}\int_{Z}s^{2}x\partial_{x}\left(  |v|^{2}\right) \\
& =\displaystyle\frac{4\tau^{2\gamma/3+1}}{\nu^{2}}\int_{Z}s^{2}%
|v|^{2}-\displaystyle\frac{4\tau^{2\gamma/3+1}}{\nu^{2}}\int_{-S_{0}}^{S_{0}%
}\left[  s^{2}x|v|^{2}\right]  _{x=0}^{x=1}\text{ ,}%
\end{array}
\]
and
\[
{\mathcal{I}}_{22}=\frac{4\tau^{2\gamma/3+1}}{\nu^{2}}\int_{Z}s^{2}%
|v|^{2}\text{ .}%
\]
Summing the above quantities yields the result, by remarking that all the
volumic terms cancel. We second compute the scalar product $\left(
{\mathcal{S}}_{x}v,{\mathcal{A}}_{s}v\right)  _{Z}$. We shall denote by
${\mathcal{J}}_{ij}$ the scalar product between the $i^{th}$ term of
${\mathcal{S}}_{x}$ with the $j^{th}$ term of ${\mathcal{A}}_{s}$.
Integrations by parts then give%
\[%
\begin{array}
[c]{ll}%
{\mathcal{J}}_{11} & =\displaystyle\frac{4\tau^{\gamma/3}}{\nu}\int
_{Z}s\partial_{x}\left(  x^{\alpha}\partial_{x}v\right)  \partial
_{s}v=-\displaystyle\frac{2\tau^{\gamma/3}}{\nu}\int_{Z}x^{\alpha}%
s\partial_{s}\left(  |\partial_{x}v|^{2}\right)  +\displaystyle\frac
{4\tau^{\gamma/3}}{\nu}\int_{-S_{0}}^{S_{0}}\left[  x^{\alpha}s\partial
_{x}v\partial_{s}v\right]  _{x=0}^{x=1}\\
& =\displaystyle\frac{2\tau^{\gamma/3}}{\nu}\int_{Z}x^{\alpha}|\partial
_{x}v|^{2}-\displaystyle\frac{2\tau^{\gamma/3}}{\nu}\int_{0}^{1}\left[
x^{\alpha}s|\partial_{x}v|^{2}\right]  _{s=-S_{0}}^{s=S_{0}}%
+\displaystyle\frac{4\tau^{\gamma/3}}{\nu}\int_{-S_{0}}^{S_{0}}\left[
x^{\alpha}s\partial_{x}v\partial_{s}v\right]  _{x=0}^{x=1}\text{ ,}%
\end{array}
\]%
\[
{\mathcal{J}}_{12}=\frac{2\tau^{\gamma/3}}{\nu}\int_{Z}v\partial_{x}\left(
x^{\alpha}\partial_{x}v\right)  =-\frac{2\tau^{\gamma/3}}{\nu}\int
_{Z}x^{\alpha}|\partial_{x}v|^{2}+\frac{2\tau^{\gamma/3}}{\nu}\int_{-S_{0}%
}^{S_{0}}\left[  x^{\alpha}v\partial_{x}v\right]  _{x=0}^{x=1}\text{ ,}%
\]%
\[%
\begin{array}
[c]{ll}%
{\mathcal{J}}_{21} & =\displaystyle\frac{4\tau^{2+\gamma/3}}{\nu}\int
_{Z}x^{2-\alpha}vs\partial_{s}v=\displaystyle\frac{2\tau^{2+\gamma/3}}{\nu
}\int_{Z}x^{2-\alpha}s\partial_{s}\left(  |v|^{2}\right) \\
& =-\displaystyle\frac{2\tau^{2+\gamma/3}}{\nu}\int_{Z}x^{2-\alpha}%
|v|^{2}+\displaystyle\frac{2\tau^{2+\gamma/3}}{\nu}\int_{0}^{1}\left[
x^{2-\alpha}s|v|^{2}\right]  _{s=-S_{0}}^{s=S_{0}}\text{ ,}%
\end{array}
\]
and%
\[
{\mathcal{J}}_{22}=\frac{2\tau^{2+\gamma/3}}{\nu}\int_{Z}x^{2-\alpha}%
|v|^{2}\text{ .}%
\]
It remains to sum the above ${\mathcal{J}}_{ij}$ to obtain the sought result
of Lemma 2.5.\hfill$\Box$

\bigskip

\section{Applications of spectral inequality}

\bigskip

The second part of this article is devoted to show some applications of the
spectral inequality.

\bigskip

Let $H$ be a real Hilbert space, and $P$ a linear self-adjoint operator from
$D(P)$ into $H$, where $D(P)$ being the domain of $P$ is a subspace of $H$.
Denote by $\left\Vert \cdot\right\Vert $ and $\left\langle \cdot
,\cdot\right\rangle $ the norm and the inner product of $H$ respectively. We
assume that $P$ is an isomorphism from $D(P)$ (equipped with the graph norm)
onto $H$, that $P^{-1}$ is a linear compact operator in $H$ and that
$\left\langle P\vartheta,\vartheta\right\rangle >0$ $\forall\vartheta\in
D(P)$, $\vartheta\neq0$. Introduce the set $\left\{  \lambda_{j}\right\}
_{j\geq1}$ for the family of all eigenvalues of $P$ so that%
\[
0<\lambda_{1}\leq\lambda_{2}\leq\cdot\cdot\leq\lambda_{k}\leq\lambda_{k+1}%
\leq\cdot\cdot\cdot\text{ and }\underset{j\rightarrow\infty}{\text{lim}%
}\lambda_{j}=\infty\text{ ,}%
\]
and let $\left\{  \Phi_{j}\right\}  _{j\geq1}$ be the family of the
corresponding orthogonal normalized eigenfunctions.

\bigskip

It is well known that for $u_{0}\in H$ given, the initial value problem
\[
\left\{
\begin{array}
[c]{ll}%
u^{\prime}\left(  t\right)  +Pu\left(  t\right)  =0\text{ ,} & t\in\left(
0,+\infty\right)  \ \text{,}\\
u\left(  0\right)  =u_{0}\text{ ,} &
\end{array}
\right.
\]
possesses a unique solution $u\in L^{2}\left(  0,T;D\left(  P^{1/2}\right)
\right)  \cap C\left(  \left[  0,T\right]  ,H\right)  $ for any $T>0$ which
satisfies%
\[
u\left(  t\right)  =\sum_{j\geq1}\left\langle u_{0},\Phi_{j}\right\rangle
e^{-\lambda_{j}t}\Phi_{j}\text{ and }\left\Vert u\left(  t\right)  \right\Vert
\leq e^{-\lambda_{1}t}\left\Vert u_{0}\right\Vert \text{ .}%
\]
In particular, if $u_{0}=\sum\limits_{j\geq1}a_{j}\Phi_{j}$ with
$\sum\limits_{j\geq1}\left\vert a_{j}\right\vert ^{2}<+\infty$, then
$\left\Vert u_{0}\right\Vert ^{2}=\sum\limits_{j\geq1}\left\vert
a_{j}\right\vert ^{2}$, $\left\langle Pu_{0},u_{0}\right\rangle =\sum
\limits_{j\geq1}\lambda_{j}\left\vert a_{j}\right\vert ^{2}$ and $\left\langle
P^{-1}u_{0},u_{0}\right\rangle =\sum\limits_{j\geq1}\frac{1}{\lambda_{j}%
}\left\vert a_{j}\right\vert ^{2}$. Further, $\frac{d}{dt}\left\Vert u\left(
t\right)  \right\Vert ^{2}+2\left\langle Pu\left(  t\right)  ,u\left(
t\right)  \right\rangle =0$ and $\frac{d}{dt}\left\langle P^{-1}u\left(
t\right)  ,u\left(  t\right)  \right\rangle +2\left\Vert u\left(  t\right)
\right\Vert ^{2}=0$.

\bigskip

Let $\Omega$ be a bounded domain of $\mathbb{R}^{d}$, $d\geq1$, with boundary
$\partial\Omega$ of class $C^{2}$. Four examples of operator $P$ are the following:

\begin{itemize}
\item The $1d$ degenerated operator with $d=1$ and $P=-\partial_{x}\left(
x^{\alpha}\partial_{x}\right)  $ with $\Omega=\left(  0,1\right)  $,
$H=L^{2}\left(  \Omega\right)  $ and $D(P)=\left\{  \vartheta\in H_{\alpha
}^{1}\left(  \Omega\right)  \text{; }{\mathcal{P}}\vartheta\in L^{2}%
(\Omega)\text{ and BC}_{\alpha}(\vartheta)=0\right\}  $ ;

\item The Laplacian with $P=-\Delta$ with $H=L^{2}\left(  \Omega\right)  $ and
$D(P)=H^{2}(\Omega)\cap H_{0}^{1}(\Omega)$ ;

\item The bi-Laplacian with $P=\Delta^{2}$ with $H=L^{2}\left(  \Omega\right)
$ and $D(P)=H^{4}(\Omega)\cap H_{0}^{2}(\Omega)$ ;

\item The Stokes operator with $P=-\mathbb{P}\Delta$ with $H=\left\{
\vartheta\in L^{2}\left(  \Omega\right)  ^{d};\text{div}\vartheta=0\text{,
}\vartheta\cdot n_{\left\vert \partial\Omega\right.  }=0\right\}  $ and
$D(P)=H^{2}(\Omega)^{d}\cap\left\{  \vartheta\in H_{0}^{1}\left(
\Omega\right)  ^{d};\text{div}\vartheta=0\right\}  $ where $\mathbb{P}$ is the
orthogonal projector in $L^{2}\left(  \Omega\right)  ^{d}$ onto $H$.
\end{itemize}

\bigskip

\subsection{Equivalence between observation and spectral inequality}

\bigskip

In this section, we present several equivalent inequalities. From now, suppose
that $H=L^{2}\left(  \Omega\right)  $. Denote $\left\Vert \cdot\right\Vert
_{\omega}$ and $\left\langle \cdot,\cdot\right\rangle _{\omega}$ the norm and
the inner product of $L^{2}\left(  \omega\right)  $ respectively where
$\omega$ is a subdomain of $\Omega$.

\bigskip

Theorem 3.1 .- \textit{Let }$\omega$\textit{ be an open and nonempty subset of
}$\Omega$\textit{. Let }$\sigma\in\left(  0,1\right)  $\textit{. Then the
following statements are equivalent:}

\begin{description}
\item[$\left(  i\right)  $] \textit{There is a positive constant }$C_{1}%
$\textit{, depending only on }$P$, $\Omega$, $\omega$\textit{ and }$\sigma
$\textit{, so that for each }$\Lambda>0$\textit{ and each sequence of real
numbers }$\left\{  a_{j}\right\}  \subset\mathbb{R}$\textit{, it holds }%
\[
\sum_{\lambda_{j}\leq\Lambda}\left\vert a_{j}\right\vert ^{2}\leq
e^{C_{1}\left(  1+\Lambda^{\sigma}\right)  }\int_{\omega}\left\vert
\sum_{\lambda_{j}\leq\Lambda}a_{j}\Phi_{j}\right\vert ^{2}\text{ .}%
\]

\item[$\left(  ii\right)  $] \textit{There is a positive constant }$C_{2}%
$\textit{, depending only on }$\left(  P,\Omega,\omega,\sigma\right)
$\textit{, so that for all }$\theta\in\left(  0,1\right)  $\textit{, }%
$t>0$\textit{ and }$u\left(  0\right)  \in L^{2}\left(  \Omega\right)
$\textit{,}%
\[
\left\Vert u\left(  t\right)  \right\Vert \leq e^{C_{2}\left(  1+\left(
\frac{1}{\theta t}\right)  ^{\frac{\sigma}{1-\sigma}}\right)  }\left\Vert
u\left(  0\right)  \right\Vert ^{\theta}\left\Vert u\left(  t\right)
\right\Vert _{\omega}^{1-\theta}\text{ .}%
\]

\item[$\left(  iii\right)  $] \textit{There is a positive constant }$C_{3}%
$\textit{, depending only on }$\left(  P,\Omega,\omega,\sigma\right)
$\textit{, so that for all }$\varepsilon>0$\textit{, }$t>0$\textit{ and
}$u\left(  0\right)  \in L^{2}\left(  \Omega\right)  $\textit{,}%
\[
\left\Vert u\left(  t\right)  \right\Vert ^{2}\leq p_{\sigma}\left(
t,\varepsilon\right)  \left\Vert u\left(  t\right)  \right\Vert _{\omega}%
^{2}+\varepsilon\left\Vert u\left(  0\right)  \right\Vert ^{2}\text{ ,}%
\]
\textit{where }%
\[
p_{\sigma}\left(  t,\varepsilon\right)  =e^{C_{3}\left(  1+\left(  \frac{1}%
{t}\right)  ^{\frac{\sigma}{1-\sigma}}\right)  }e^{\left(  \frac{C_{3}}%
{t}\text{ln}\left(  e+\frac{1}{\varepsilon}\right)  \right)  ^{\sigma}}\text{
.}%
\]

\item[$\left(  iv\right)  $] \textit{There is a positive constant }$C_{4}%
$\textit{, depending only on }$\left(  P,\Omega,\omega,\sigma\right)
$\textit{, so that for all }$t>0$\textit{ and }$u\left(  0\right)  \in
L^{2}\left(  \Omega\right)  $\textit{,}%
\[
\left\Vert u\left(  t\right)  \right\Vert \leq e^{C_{4}\left(  1+\left(
\frac{1}{t}\right)  ^{\frac{\sigma}{1-\sigma}}\right)  }e^{\left(  \frac
{C_{4}}{t}\text{ln}\left(  \frac{\left\Vert u\left(  0\right)  \right\Vert
}{\left\Vert u\left(  t\right)  \right\Vert }\right)  \right)  ^{\sigma}%
}\left\Vert u\left(  t\right)  \right\Vert _{\omega}\text{ .}%
\]

\end{description}

\bigskip

In particular, if $P=-\Delta$, then $\sigma=\frac{1}{2}$ (see \cite{L},
\cite{LZ}, \cite{PWX}, \cite{BP}, \cite{Ph2}); If $P=\Delta^{2}$, then
$\sigma=\frac{1}{4}$ (see \cite{AE}, \cite{EMZ}, \cite{Ga}, \cite{LRR2}); If
$P$ is the Stokes operator, then $\sigma=\frac{1}{2}$ (see \cite{CSL}).

\bigskip

Proof .- We organize the proof by several steps.

\noindent\textit{Step 1: To show that }$(i)\Rightarrow(ii)$.

Arbitrarily fix $\lambda>0$, $t>0$ and $u\left(  0\right)  =\sum
\limits_{j\geq1}a_{j}e_{j}$ with $\{a_{j}\}_{j\geq1}\subset\ell^{2}$. Write
\[
u\left(  t\right)  =\sum_{\lambda_{j}\leq\Lambda}a_{j}e^{-\lambda_{j}t}%
\Phi_{j}+\sum_{\lambda_{j}>\Lambda}a_{j}e^{-\lambda_{j}t}\Phi_{j}\text{ .}%
\]
Then by $(i)$, we find that%
\[%
\begin{array}
[c]{ll}%
\left\Vert u\left(  t\right)  \right\Vert  & \leq\left\Vert \displaystyle\sum
_{\lambda_{j}\leq\Lambda}a_{j}e^{-\lambda_{j}t}\Phi_{j}\right\Vert +\left\Vert
\displaystyle\sum_{\lambda_{j}>\Lambda}a_{j}e^{-\lambda_{j}t}\Phi
_{j}\right\Vert \\
& \leq\left(  \displaystyle\sum_{\lambda_{j}\leq\Lambda}\left\vert
a_{j}e^{-\lambda_{j}t}\right\vert ^{2}\right)  ^{1/2}+e^{-\Lambda t}\left\Vert
u\left(  0\right)  \right\Vert \\
& \leq\left(  e^{C_{1}\left(  1+\Lambda^{\sigma}\right)  }\displaystyle\int
_{\omega}\left\vert \displaystyle\sum_{\lambda_{j}\leq\Lambda}a_{j}%
e^{-\lambda_{j}t}\Phi_{j}\right\vert ^{2}\right)  ^{1/2}+e^{-\Lambda
t}\left\Vert u\left(  0\right)  \right\Vert \text{ .}%
\end{array}
\]
This, along with the triangle inequality for the norm $\left\Vert
\cdot\right\Vert _{\omega}$, yields that%
\[%
\begin{array}
[c]{ll}%
\left\Vert u\left(  t\right)  \right\Vert  & \leq\left(  e^{C_{1}\left(
1+\Lambda^{\sigma}\right)  }\displaystyle\int_{\omega}\left\vert
\displaystyle\sum_{j\geq1}a_{j}e^{-\lambda_{j}t}\Phi_{j}\right\vert
^{2}\right)  ^{1/2}\\
& \quad+\left(  e^{C_{1}\left(  1+\Lambda^{\sigma}\right)  }\displaystyle\int
_{\omega}\left\vert \displaystyle\sum_{\lambda_{j}>\Lambda}a_{j}%
e^{-\lambda_{j}t}\Phi_{j}\right\vert ^{2}\right)  ^{1/2}+e^{-\Lambda
t}\left\Vert u\left(  0\right)  \right\Vert \text{ .}%
\end{array}
\]
Hence, it follows that%
\[%
\begin{array}
[c]{ll}%
\left\Vert u\left(  t\right)  \right\Vert  & \leq e^{\frac{C_{1}}{2}\left(
1+\Lambda^{\sigma}\right)  }\left\Vert u\left(  t\right)  \right\Vert
_{\omega}+e^{\frac{C_{1}}{2}\left(  1+\Lambda^{\sigma}\right)  }e^{-\Lambda
t}\left\Vert u\left(  0\right)  \right\Vert +e^{-\Lambda t}\left\Vert u\left(
0\right)  \right\Vert \\
& \leq2e^{\frac{C_{1}}{2}\left(  1+\Lambda^{\sigma}\right)  }\left(
\left\Vert u\left(  t\right)  \right\Vert _{\omega}+e^{-\Lambda t}\left\Vert
u\left(  0\right)  \right\Vert \right)  \text{ .}%
\end{array}
\]
Since by the Young inequality
\[
C_{1}\Lambda^{\sigma}=\frac{C_{1}}{\left(  \epsilon t\right)  ^{\sigma}%
}\left(  \epsilon\Lambda t\right)  ^{\sigma}\leq\epsilon\Lambda t+\left(
\frac{C_{1}}{\left(  \epsilon t\right)  ^{\sigma}}\right)  ^{\frac{1}%
{1-\sigma}}\quad\text{for any }\epsilon,t>0\text{ ,}%
\]
one deduce that for all $\epsilon\in(0,2)$,
\[
\left\Vert u\left(  t\right)  \right\Vert \leq2e^{\frac{C_{1}}{2}}e^{\frac
{1}{2}\left(  \frac{C_{1}}{\left(  \epsilon t\right)  ^{\sigma}}\right)
^{\frac{1}{1-\sigma}}}\left(  e^{\frac{\epsilon}{2}\Lambda t}\left\Vert
u\left(  t\right)  \right\Vert _{\omega}+e^{-\frac{2-\epsilon}{2}\Lambda
t}\left\Vert u\left(  0\right)  \right\Vert \right)  \quad\text{\ for each
}\Lambda>0\text{ .}%
\]
Notice that if $\left\Vert u\left(  t\right)  \right\Vert _{\omega}=0$ then,
$\left\Vert u\left(  t\right)  \right\Vert =0$. Next, choose
\[
\Lambda=\frac{1}{t}\text{ln}\left(  \frac{\left\Vert u\left(  0\right)
\right\Vert }{\left\Vert u\left(  t\right)  \right\Vert _{\omega}}\right)
\]
(knowing that $\left\Vert u\left(  t\right)  \right\Vert _{\omega}%
\leq\left\Vert u\left(  0\right)  \right\Vert $) to get
\[
\left\Vert u\left(  t\right)  \right\Vert \leq2e^{\frac{C_{1}}{2}}e^{\frac
{1}{2}\left(  \frac{C_{1}}{\left(  \epsilon t\right)  ^{\sigma}}\right)
^{\frac{1}{1-\sigma}}}\left(  2\left\Vert u\left(  t\right)  \right\Vert
_{\omega}^{1-\frac{\epsilon}{2}}\left\Vert u\left(  0\right)  \right\Vert
^{\frac{\epsilon}{2}}\right)
\]
which is the inequality in $(ii)$ with $\theta=\frac{\epsilon}{2}$ and
ln$4+\frac{C_{1}}{2}+\frac{1}{2}\left(  \frac{C_{1}}{\left(  \varepsilon
t\right)  ^{\sigma}}\right)  ^{\frac{1}{1-\sigma}}\leq C_{2}\left(  1+\left(
\frac{1}{\theta t}\right)  ^{\frac{\sigma}{1-\sigma}}\right)  $.

\noindent\textit{Step 2: To show that }$(ii)\Rightarrow(iii)$.

We write the inequality in $(ii)$ in the following way
\[
\left\Vert u\left(  t\right)  \right\Vert ^{2}\leq\left\Vert u\left(
0\right)  \right\Vert ^{2\theta}\left(  \text{exp}\left(  \frac{2C_{2}%
}{1-\theta}\left(  1+\left(  \frac{1}{\theta t}\right)  ^{\frac{\sigma
}{1-\sigma}}\right)  \right)  \left\Vert u\left(  t\right)  \right\Vert
_{\omega}^{2}\right)  ^{1-\theta}%
\]
and apply the fact that for any $E,B,D>0$ and $\theta\in(0,1)$
\[
E\leq B^{\theta}D^{1-\theta}\Leftrightarrow E\leq\varepsilon B+\left(
1-\theta\right)  \theta^{\frac{\theta}{1-\theta}}\frac{1}{\varepsilon
^{\frac{\theta}{1-\theta}}}D\quad\forall\varepsilon>0\text{ .}%
\]
To prove the above equivalence, one uses the Young inequality and one choose
$\varepsilon=\theta\left(  \frac{D}{B}\right)  ^{1-\theta}$. Therefore,
\[
\left\Vert u\left(  t\right)  \right\Vert ^{2}\leq\varepsilon\left\Vert
u\left(  0\right)  \right\Vert ^{2}+\text{exp}\left(  \frac{2C_{2}}{1-\theta
}\left(  1+\left(  \frac{1}{\theta t}\right)  ^{\frac{\sigma}{1-\sigma}%
}\right)  \right)  \frac{1}{\varepsilon^{\frac{\theta}{1-\theta}}}\left\Vert
u\left(  t\right)  \right\Vert _{\omega}^{2}\text{ .}%
\]
By denoting $\beta=\frac{\theta}{1-\theta}$, it yields%
\[
\left\Vert u\left(  t\right)  \right\Vert ^{2}\leq\varepsilon\left\Vert
u\left(  0\right)  \right\Vert ^{2}+e^{2C_{2}\left(  1+\beta\right)  \left(
1+\left(  \frac{1+\beta}{\beta t}\right)  ^{\frac{\sigma}{1-\sigma}}\right)
}\frac{1}{\varepsilon^{\beta}}\left\Vert u\left(  t\right)  \right\Vert
_{\omega}^{2}\text{ .}%
\]
Now, notice, with $B=K\left(  1+\left(  \frac{1}{t}\right)  ^{\frac{\sigma
}{1-\sigma}}\right)  $ and $D=K\left(  \frac{1}{t}\right)  ^{\frac{\sigma
}{1-\sigma}}$ for some constant $K>0$, that
\[
\frac{1}{\varepsilon^{\beta}}e^{2C_{2}\left(  1+\beta\right)  \left(
1+\left(  \frac{1+\beta}{\beta t}\right)  ^{\frac{\sigma}{1-\sigma}}\right)
}\leq e^{B+\beta\left(  \text{ln}\left(  e+\frac{1}{\varepsilon}\right)
+B\right)  +\left(  \frac{1}{\beta}\right)  ^{\frac{\sigma}{1-\sigma}}D}\text{
.}%
\]
Next, choose $\beta=\left(  \frac{D}{\text{ln}\left(  e+\frac{1}{\varepsilon
}\right)  +B}\right)  ^{1-\sigma}$ to get
\[
e^{B+\beta\left(  \text{ln}\left(  e+\frac{1}{\varepsilon}\right)  +B\right)
+\left(  \frac{1}{\beta}\right)  ^{\frac{\sigma}{1-\sigma}}D}\leq
e^{cB+c\left(  \text{ln}\left(  e+\frac{1}{\varepsilon}\right)  +B\right)
^{\sigma}D^{1-\sigma}}\leq e^{c^{\prime}B+c^{\prime}\left(  \text{ln}\left(
e+\frac{1}{\varepsilon}\right)  \right)  ^{\sigma}D^{1-\sigma}}%
\]
for some constants $c,c^{\prime}>0$. Therefore, we obtain the desired
inequality
\[
\left\Vert u\left(  t\right)  \right\Vert ^{2}\leq\varepsilon\left\Vert
u\left(  0\right)  \right\Vert ^{2}+e^{C_{3}\left(  1+\left(  \frac{1}%
{t}\right)  ^{\frac{\sigma}{1-\sigma}}\right)  }e^{\left(  \frac{C_{3}}%
{t}\text{ln}\left(  e+\frac{1}{\varepsilon}\right)  \right)  ^{\sigma}%
}\left\Vert u\left(  t\right)  \right\Vert _{\omega}^{2}\text{ .}%
\]

\noindent\textit{Step 3: to show that }$(iii)\Rightarrow(iv)$.

Take
\[
\varepsilon=\frac{1}{2}\frac{\left\Vert u\left(  t\right)  \right\Vert ^{2}%
}{\left\Vert u\left(  0\right)  \right\Vert ^{2}}%
\]
in the inequality in $(iii)$ and we use the fact that $\left\Vert u\left(
t\right)  \right\Vert \leq\left\Vert u\left(  0\right)  \right\Vert $.
Therefore, we have
\[
\frac{1}{2}\left\Vert u\left(  t\right)  \right\Vert ^{2}\leq e^{C_{3}\left(
1+\left(  \frac{1}{t}\right)  ^{\frac{\sigma}{1-\sigma}}\right)  }e^{\left(
\frac{C_{3}}{t}\text{ln}\left(  \left(  e+2\right)  \frac{\left\Vert u\left(
0\right)  \right\Vert ^{2}}{\left\Vert u\left(  t\right)  \right\Vert ^{2}%
}\right)  \right)  ^{\sigma}}\left\Vert u\left(  t\right)  \right\Vert
_{\omega}^{2}\text{ .}%
\]

\noindent\textit{Step 4: to show that }$(iv)\Rightarrow(i)$.

Apply the Young inequality%
\[
\left(  \frac{C_{4}}{t}\text{ln}\left(  \frac{\left\Vert u\left(  0\right)
\right\Vert }{\left\Vert u\left(  t\right)  \right\Vert }\right)  \right)
^{\sigma}\leq\left(  \frac{C_{4}}{t}\right)  ^{\frac{\sigma}{1-\sigma}%
}+\text{ln}\left(  \frac{\left\Vert u\left(  0\right)  \right\Vert
}{\left\Vert u\left(  t\right)  \right\Vert }\right)
\]
to deduce the inequality: There are two constants $C>0$ and $\alpha\in\left(
0,1\right)  $, which depend only on $\left(  \Omega,\omega,\sigma\right)  $,
so that for all $t>0$ and $u\left(  0\right)  \in L^{2}\left(  \Omega\right)
$,%
\[
\left\Vert u\left(  t\right)  \right\Vert \leq e^{C\left(  1+\left(  \frac
{1}{t}\right)  ^{\frac{\sigma}{1-\sigma}}\right)  }\left\Vert u\left(
0\right)  \right\Vert ^{\alpha}\left\Vert u\left(  t\right)  \right\Vert
_{\omega}^{1-\alpha}\text{ .}%
\]
Arbitrarily fix $\Lambda>0$ and $\{a_{j}\}\subset\mathbb{R}$. By applying the
above inequality, with $u\left(  0\right)  =\sum\limits_{\lambda_{j}<\Lambda
}a_{j}e^{\lambda_{j}t}\Phi_{j}$, we get that
\[
\sum_{\lambda_{j}\leq\Lambda}\left\vert a_{j}\right\vert ^{2}\leq e^{2C\left(
1+\left(  \frac{1}{t}\right)  ^{\frac{\sigma}{1-\sigma}}\right)  }\left(
\sum_{\lambda_{j}\leq\Lambda}\left\vert a_{j}e^{\lambda_{j}t}\right\vert
^{2}\right)  ^{\alpha}\left(  \int_{\omega}\left\vert \sum_{\lambda_{j}%
\leq\Lambda}a_{j}\Phi_{j}\right\vert ^{2}\right)  ^{1-\alpha}\text{ ,}%
\]
which implies that
\[
\sum_{\lambda_{j}\leq\Lambda}|a_{j}|^{2}\leq e^{\frac{2}{1-\alpha}C\left(
1+\left(  \frac{1}{t}\right)  ^{\frac{\sigma}{1-\sigma}}\right)  }%
e^{\frac{2\alpha}{1-\alpha}\Lambda t}\int_{\omega}\left\vert \sum_{\lambda
_{j}\leq\Lambda}a_{j}\Phi_{j}\right\vert ^{2}\quad\text{\ for each }t>0\text{
.}%
\]
Choose $t=\left(  \frac{1}{\Lambda}\right)  ^{1-\sigma}$ to get the conclusion
$(i)$.

\bigskip

This ends the proof.\hfill$\Box$

\bigskip

\subsection{Equivalence between observation and control}

\bigskip

Let us recall the classical results of equivalence between observation
estimate and controllability with cost. There are at least three ways to
establish the cost: One is based on the duality of the control operator in the
spirit of the HUM method (see \cite{Lio}) with a spectral decomposition (see
\cite{R}, \cite{Ph}); Another one have a geometric point of view using
Hahn-Banach Theorem (see \cite{WWZ}, \cite{WYZ}) ; The last one is based on a
minimization of a certain functional (see \cite{FZ}, \cite{Mi}). The arguments
we present are similar to those appear in \cite[lemma 3.2, p.1475]{Mi} (see
also \cite[remark 6.6, p.3670]{DM}).

\bigskip

Denote $\left\Vert \cdot\right\Vert $ and $\left\langle \cdot,\cdot
\right\rangle $ the norm and the inner product of $L^{2}\left(  \Omega\right)
$ respectively.

\bigskip

Theorem 3.2.- \textit{Let }$0\leq T_{0}<T_{1}<T_{2}$. \textit{Let }%
$\ell,\varepsilon>0$\textit{. The following two statements are equivalent.}

\begin{description}
\item[\textit{(}$\mathcal{C}$\textit{)}] \textit{For any }$y_{e}\in
L^{2}\left(  \Omega\right)  $\textit{, there is }$f\in L^{2}\left(
\omega\right)  $\textit{ such that the solution} $y$\textit{ to}
\[
\left\{
\begin{array}
[c]{ll}%
y^{\prime}\left(  t\right)  +Py\left(  t\right)  =0\text{ ,} & t\in\left(
T_{0},T_{2}\right)  \backslash\left\{  T_{1}\right\}  \ \text{,}\\
y\left(  T_{0}\right)  =y_{e}\text{ ,} & \\
y\left(  T_{1}\right)  =y\left(  T_{1-}\right)  +1_{\omega}f\text{ ,} &
\end{array}
\right.
\]
\textit{satisfies }%
\[
\frac{1}{\ell}\left\Vert f\right\Vert _{\omega}^{2}+\frac{1}{\varepsilon
}\left\Vert y\left(  T_{2}\right)  \right\Vert ^{2}\leq\left\Vert
y_{e}\right\Vert ^{2}\text{ .}%
\]

\item[\textit{(}$\mathcal{O}$\textit{)}] \textit{The solution }$u$\textit{
to}
\[
\left\{
\begin{array}
[c]{ll}%
u^{\prime}\left(  t\right)  +Pu\left(  t\right)  =0\text{ ,} & t\in\left(
T_{0},T_{2}\right)  \ \text{,}\\
u\left(  T_{0}\right)  \in L^{2}\left(  \Omega\right)  \text{ ,} &
\end{array}
\right.
\]
\textit{satisfies}%
\[
\left\Vert u\left(  T_{2}\right)  \right\Vert ^{2}\leq\ell\left\Vert u\left(
T_{0}+T_{2}-T_{1}\right)  \right\Vert _{\omega}^{2}+\varepsilon\left\Vert
u\left(  T_{0}\right)  \right\Vert ^{2}\text{ .}%
\]

\end{description}

\bigskip

Proof of $(\mathcal{C})\Rightarrow(\mathcal{O})$ .- We multiply the equations
of $(\mathcal{C})$ by $u\left(  T_{0}+T_{2}-t\right)  $ to get
\[
\left\langle y(T_{2}),u\left(  T_{0}\right)  \right\rangle -\left\langle
y(T_{0}),u\left(  T_{2}\right)  \right\rangle =\left\langle f,u\left(
T_{0}+T_{2}-T_{1}\right)  \right\rangle _{\omega}\text{ ,}%
\]
that is,
\[
\left\langle y_{e},u\left(  T_{2}\right)  \right\rangle =-\left\langle
f,u\left(  T_{0}+T_{2}-T_{1}\right)  \right\rangle _{\omega}+\left\langle
y(T_{2}),u\left(  T_{0}\right)  \right\rangle \text{ .}%
\]
By Cauchy-Schwarz inequality and using the inequality in $(\mathcal{C})$ one
can deduce that
\[%
\begin{array}
[c]{ll}%
\left\langle y_{e},u\left(  T_{2}\right)  \right\rangle  & \leq\left\Vert
f\right\Vert _{\omega}\left\Vert u\left(  T_{0}+T_{2}-T_{1}\right)
\right\Vert _{\omega}+\left\Vert y\left(  T_{2}\right)  \right\Vert \left\Vert
u\left(  T_{0}\right)  \right\Vert \\
& \leq\dfrac{1}{2\ell}\left\Vert f\right\Vert _{\omega}^{2}+\dfrac
{1}{2\varepsilon}\left\Vert y\left(  T_{2}\right)  \right\Vert ^{2}%
+\dfrac{\ell}{2}\left\Vert u\left(  T_{0}+T_{2}-T_{1}\right)  \right\Vert
_{\omega}^{2}+\dfrac{\varepsilon}{2}\left\Vert u\left(  T_{0}\right)
\right\Vert ^{2}\\
& \leq\dfrac{1}{2}\left\Vert y_{e}\right\Vert ^{2}+\dfrac{1}{2}\left(
\ell\left\Vert u\left(  T_{0}+T_{2}-T_{1}\right)  \right\Vert _{\omega}%
^{2}+\varepsilon\left\Vert u\left(  T_{0}\right)  \right\Vert ^{2}\right)
\end{array}
\]
which gives the desired estimate by choosing $y_{e}=u\left(  T_{2}\right)  $.

\bigskip

Proof of $(\mathcal{O})\Rightarrow(\mathcal{C})$ .- Let $y_{e}\in L^{2}\left(
\Omega\right)  $. Consider the functional $J$ defined on $L^{2}\left(
\Omega\right)  $ given by%
\[
J\left(  \vartheta\right)  =\frac{\ell}{2}\left\Vert u\left(  T_{0}%
+T_{2}-T_{1}\right)  \right\Vert _{\omega}^{2}+\frac{\varepsilon}{2}\left\Vert
\vartheta\right\Vert ^{2}-\left\langle y_{e},u\left(  T_{2}\right)
\right\rangle \text{ ,}%
\]
where%
\[
\left\{
\begin{array}
[c]{ll}%
u^{\prime}\left(  t\right)  +Pu\left(  t\right)  =0\text{ ,} & t\in\left(
T_{0},T_{2}\right)  \ \text{,}\\
u\left(  T_{0}\right)  =\vartheta\text{ .} &
\end{array}
\right.
\]
Notice that $J$ is strictly convex, $C^{1}$ and coercive and therefore $J$ has
a unique minimizer $w_{0}\in L^{2}\left(  \Omega\right)  $, i.e.
$J(w_{0})=\underset{\vartheta\in L^{2}\left(  \Omega\right)  }{\text{min}%
}J(\vartheta)$. Set%
\[
\left\{
\begin{array}
[c]{ll}%
w^{\prime}\left(  t\right)  +Pw\left(  t\right)  =0\text{ ,} & t\in\left(
T_{0},T_{2}\right)  \ \text{,}\\
w\left(  T_{0}\right)  =w_{0}\text{ ,} &
\end{array}
\right.  \text{ and }\left\{
\begin{array}
[c]{ll}%
h^{\prime}\left(  t\right)  +Ph\left(  t\right)  =0\text{ ,} & t\in\left(
T_{0},T_{2}\right)  \ \text{,}\\
h\left(  T_{0}\right)  =h_{0}\text{ .} &
\end{array}
\right.
\]
Since $J^{\prime}(w_{0})h_{0}=0$ for any $h_{0}\in L^{2}\left(  \Omega\right)
$, we have%
\[
\ell\left\langle w\left(  T_{0}+T_{2}-T_{1}\right)  ,h\left(  T_{0}%
+T_{2}-T_{1}\right)  \right\rangle _{\omega}+\varepsilon\left\langle
w_{0},h_{0}\right\rangle -\left\langle y_{e},h\left(  T_{2}\right)
\right\rangle =0\quad\forall h_{0}\in L^{2}\left(  \Omega\right)  \text{ .}%
\]
On the other hand, the identity
\[
\left\langle y\left(  T_{2}\right)  ,u\left(  T_{0}\right)  \right\rangle
-\left\langle y_{e},u\left(  T_{2}\right)  \right\rangle =\left\langle
f,u\left(  T_{0}+T_{2}-T_{1}\right)  \right\rangle _{\omega}\quad\forall
u\left(  T_{0}\right)  \in L^{2}\left(  \Omega\right)
\]
implies
\[
-\left\langle f,h\left(  T_{0}+T_{2}-T_{1}\right)  \right\rangle _{\omega
}+\left\langle y\left(  T_{2}\right)  ,h_{0}\right\rangle -\left\langle
y_{e},h\left(  T_{2}\right)  \right\rangle =0\quad\forall h_{0}\in
L^{2}\left(  \Omega\right)  \text{ .}%
\]
By choosing $f=-\ell w\left(  T_{0}+T_{2}-T_{1}\right)  $, we deduce that the
solution $y$ satisfies
\[
\varepsilon w_{0}=y\left(  T_{2}\right)  \text{ .}%
\]
Further,
\[
\ell\left\Vert w\left(  T_{0}+T_{2}-T_{1}\right)  \right\Vert _{\omega}%
^{2}+\varepsilon\left\Vert w_{0}\right\Vert ^{2}=\frac{1}{\ell}\left\Vert
f\right\Vert _{\omega}^{2}+\frac{1}{\varepsilon}\left\Vert y\left(
T_{2}\right)  \right\Vert ^{2}\text{ .}%
\]
Moreover, taking $h_{0}=w_{0}$ into $J^{\prime}(w_{0})h_{0}=0$, we get%
\[
\ell\left\Vert w\left(  T_{0}+T_{2}-T_{1}\right)  \right\Vert _{\omega}%
^{2}+\varepsilon\left\Vert w_{0}\right\Vert ^{2}-\left\langle y_{e},w\left(
T_{2}\right)  \right\rangle =0\text{ .}%
\]
By Cauchy-Schwarz inequality,
\[%
\begin{array}
[c]{ll}%
\ell\left\Vert w\left(  T_{0}+T_{2}-T_{1}\right)  \right\Vert _{\omega}%
^{2}+\varepsilon\left\Vert w_{0}\right\Vert ^{2} & \leq\left\Vert
y_{e}\right\Vert _{L^{2}\left(  \Omega\right)  }\left\Vert w\left(
T_{2}\right)  \right\Vert _{L^{2}\left(  \Omega\right)  }\\
& \leq\left\Vert y_{e}\right\Vert _{L^{2}\left(  \Omega\right)  }\left(
\ell\left\Vert w\left(  T_{0}+T_{2}-T_{1}\right)  \right\Vert _{\omega}%
^{2}+\varepsilon\left\Vert w_{0}\right\Vert ^{2}\right)  ^{1/2}%
\end{array}
\]
where in the last line, we used $(\mathcal{O})$. Therefore, we get%
\[
\ell\left\Vert w\left(  T_{0}+T_{2}-T_{1}\right)  \right\Vert _{\omega}%
^{2}+\varepsilon\left\Vert w_{0}\right\Vert ^{2}\leq\left\Vert y_{e}%
\right\Vert ^{2}\text{ ,}%
\]
that is,
\[
\frac{1}{\ell}\left\Vert f\right\Vert _{\omega}^{2}+\frac{1}{\varepsilon
}\left\Vert y\left(  T_{2}\right)  \right\Vert ^{2}\leq\left\Vert
y_{e}\right\Vert ^{2}\text{ }%
\]
where%
\[
\left\{
\begin{array}
[c]{ll}%
y^{\prime}\left(  t\right)  +Py\left(  t\right)  =0\text{ ,} & t\in\left(
T_{0},T_{2}\right)  \backslash\left\{  T_{1}\right\}  \ \text{,}\\
y\left(  T_{0}\right)  =y_{e}\text{ ,} & \\
y\left(  T_{1}\right)  =y\left(  T_{1-}\right)  +1_{\omega}\left(  -\ell
w\left(  T_{0}+T_{2}-t\right)  \right)  \text{ ,} & \\
w^{\prime}\left(  t\right)  +Pw\left(  t\right)  =0\text{ ,} & t\in\left(
T_{0},T_{2}\right)  \ \text{,}\\
w\left(  T_{0}\right)  =\frac{1}{\varepsilon}y\left(  T_{2}\right)  \text{ .}
&
\end{array}
\right.
\]
This completes the proof. \hfill$\Box$

\bigskip

\bigskip

\bigskip

Theorem 3.3.- \textit{Let }$0\leq T_{0}<T_{1}<T_{2}$. \textit{Let }%
$\ell,\varepsilon>0$\textit{. The following two statements are equivalent.}

\begin{description}
\item[\textit{(}$\mathcal{C}$\textit{)}] \textit{For any }$y_{d}\in
L^{2}\left(  \Omega\right)  $\textit{ such that }$\left\langle Py_{d}%
,y_{d}\right\rangle <+\infty$\textit{, there is }$f\in L^{2}\left(
\omega\right)  $\textit{ such that the solution} $y$\textit{ to}
\[
\left\{
\begin{array}
[c]{ll}%
y^{\prime}\left(  t\right)  +Py\left(  t\right)  =0\text{ ,} & t\in\left(
T_{0},T_{2}\right)  \backslash\left\{  T_{1}\right\}  \ \text{,}\\
y\left(  T_{0}\right)  =0\text{ ,} & \\
y\left(  T_{1}\right)  =y\left(  T_{1-}\right)  +1_{\omega}f\text{ ,} &
\end{array}
\right.
\]
\textit{satisfies }%
\[
\frac{1}{\ell}\left\Vert f\right\Vert _{\omega}^{2}+\frac{1}{\varepsilon
}\left\Vert y\left(  T_{2}\right)  -y_{d}\right\Vert ^{2}\leq\left\langle
Ay_{d},y_{d}\right\rangle \text{ .}%
\]

\item[\textit{(}$\mathcal{O}$\textit{)}] \textit{The solution }$u$\textit{
to}
\[
\left\{
\begin{array}
[c]{ll}%
u^{\prime}\left(  t\right)  +Pu\left(  t\right)  =0\text{ ,} & t\in\left(
T_{0},T_{2}\right)  \ \text{,}\\
u\left(  T_{0}\right)  \in L^{2}\left(  \Omega\right)  \text{ ,} &
\end{array}
\right.
\]
\textit{satisfies}%
\[
\left\langle P^{-1}u\left(  T_{0}\right)  ,u\left(  T_{0}\right)
\right\rangle \leq\ell\left\Vert u\left(  T_{0}+T_{2}-T_{1}\right)
\right\Vert _{\omega}^{2}+\varepsilon\left\Vert u\left(  T_{0}\right)
\right\Vert ^{2}\text{ .}%
\]

\end{description}

\bigskip

Proof of $(\mathcal{C})\Rightarrow(\mathcal{O})$ .- We multiply the equations
of $(\mathcal{C})$ by $u\left(  T_{0}+T_{2}-t\right)  $ to get
\[
\left\langle y(T_{2}),u\left(  T_{0}\right)  \right\rangle -\left\langle
y(T_{0}),u\left(  T_{2}\right)  \right\rangle =\left\langle f,u\left(
T_{0}+T_{2}-T_{1}\right)  \right\rangle _{\omega}\text{ ,}%
\]
that is,
\[
\left\langle y_{d},u\left(  T_{0}\right)  \right\rangle =\left\langle
f,u\left(  T_{0}+T_{2}-T_{1}\right)  \right\rangle _{\omega}-\left\langle
y(T_{2})-y_{d},u\left(  T_{0}\right)  \right\rangle \text{ .}%
\]
By Cauchy-Schwarz inequality and using the inequality in $(\mathcal{C})$ one
has
\[%
\begin{array}
[c]{ll}%
\left\langle y_{d},u\left(  T_{0}\right)  \right\rangle  & \leq\left\Vert
f\right\Vert _{\omega}\left\Vert u\left(  T_{0}+T_{2}-T_{1}\right)
\right\Vert _{\omega}+\left\Vert y\left(  T_{2}\right)  -y_{d}\right\Vert
\left\Vert u\left(  T_{0}\right)  \right\Vert \\
& \leq\dfrac{1}{2\ell}\left\Vert f\right\Vert _{\omega}^{2}+\dfrac
{1}{2\varepsilon}\left\Vert y\left(  T_{2}\right)  -y_{d}\right\Vert
^{2}+\dfrac{\ell}{2}\left\Vert u\left(  T_{0}+T_{2}-T_{1}\right)  \right\Vert
_{\omega}^{2}+\dfrac{\varepsilon}{2}\left\Vert u\left(  T_{0}\right)
\right\Vert ^{2}\\
& \leq\dfrac{1}{2}\left\langle Py_{d},y_{d}\right\rangle +\dfrac{1}{2}\left(
\ell\left\Vert u\left(  T_{0}+T_{2}-T_{1}\right)  \right\Vert _{\omega}%
^{2}+\varepsilon\left\Vert u\left(  T_{0}\right)  \right\Vert ^{2}\right)
\end{array}
\]
which gives the desired estimate by choosing $y_{d}=P^{-1}u\left(
T_{0}\right)  $.

\bigskip

Proof of $(\mathcal{O})\Rightarrow(\mathcal{C})$ .- Let $y_{d}\in L^{2}\left(
\Omega\right)  $ such that $\left\langle Py_{d},y_{d}\right\rangle <+\infty$.
Consider the functional $J$ defined on $L^{2}\left(  \Omega\right)  $ given by%
\[
J\left(  \vartheta\right)  =\frac{\ell}{2}\left\Vert u\left(  T_{0}%
+T_{2}-T_{1}\right)  \right\Vert _{\omega}^{2}+\frac{\varepsilon}{2}\left\Vert
\vartheta\right\Vert ^{2}+\left\langle y_{d},\vartheta\right\rangle \text{ ,}%
\]
where%
\[
\left\{
\begin{array}
[c]{ll}%
u^{\prime}\left(  t\right)  +Pu\left(  t\right)  =0\text{ ,} & t\in\left(
T_{0},T_{2}\right)  \ \text{,}\\
u\left(  T_{0}\right)  =\vartheta\text{ .} &
\end{array}
\right.
\]
Notice that $J$ is strictly convex, $C^{1}$ and coercive and therefore $J$ has
a unique minimizer $w_{0}\in L^{2}\left(  \Omega\right)  $, i.e.
$J(w_{0})=\underset{\vartheta\in L^{2}\left(  \Omega\right)  }{\text{min}%
}J(\vartheta)$. Set%
\[
\left\{
\begin{array}
[c]{ll}%
w^{\prime}\left(  t\right)  +Pw\left(  t\right)  =0\text{ ,} & t\in\left(
T_{0},T_{2}\right)  \ \text{,}\\
w\left(  T_{0}\right)  =w_{0}\text{ ,} &
\end{array}
\right.  \text{ and }\left\{
\begin{array}
[c]{ll}%
h^{\prime}\left(  t\right)  +Ph\left(  t\right)  =0\text{ ,} & t\in\left(
T_{0},T_{2}\right)  \ \text{,}\\
h\left(  T_{0}\right)  =h_{0}\text{ .} &
\end{array}
\right.
\]
Since $J^{\prime}(w_{0})h_{0}=0$ for any $h_{0}\in L^{2}\left(  \Omega\right)
$, we have%
\[
\ell\left\langle w\left(  T_{0}+T_{2}-T_{1}\right)  ,h\left(  T_{0}%
+T_{2}-T_{1}\right)  \right\rangle _{\omega}+\varepsilon\left\langle
w_{0},h_{0}\right\rangle +\left\langle y_{d},h_{0}\right\rangle =0\quad\forall
h_{0}\in L^{2}\left(  \Omega\right)  \text{ .}%
\]
On the other hand, the identity
\[
\left\langle y\left(  T_{2}\right)  ,u\left(  T_{0}\right)  \right\rangle
-\left\langle y\left(  T_{0}\right)  ,u\left(  T_{2}\right)  \right\rangle
=\left\langle f,u\left(  T_{0}+T_{2}-T_{1}\right)  \right\rangle _{\omega
}\quad\forall u\left(  T_{0}\right)  \in L^{2}\left(  \Omega\right)
\]
implies
\[
-\left\langle f,h\left(  T_{0}+T_{2}-T_{1}\right)  \right\rangle _{\omega
}+\left\langle y\left(  T_{2}\right)  -y_{d},h_{0}\right\rangle +\left\langle
y_{d},h_{0}\right\rangle =0\quad\forall h_{0}\in L^{2}\left(  \Omega\right)
\text{ .}%
\]
By choosing $f=-\ell w\left(  T_{0}+T_{2}-T_{1}\right)  $, we deduce that the
solution $y$ satisfies
\[
\varepsilon w_{0}=y\left(  T_{2}\right)  -y_{d}\text{ .}%
\]
Further,
\[
\ell\left\Vert w\left(  T_{0}+T_{2}-T_{1}\right)  \right\Vert _{\omega}%
^{2}+\varepsilon\left\Vert w_{0}\right\Vert ^{2}=\frac{1}{\ell}\left\Vert
f\right\Vert _{\omega}^{2}+\frac{1}{\varepsilon}\left\Vert y\left(
T_{2}\right)  -y_{d}\right\Vert ^{2}\text{ .}%
\]
Moreover, taking $h_{0}=w_{0}$ into $J^{\prime}(w_{0})h_{0}=0$, we get%
\[
\ell\left\Vert w\left(  T_{0}+T_{2}-T_{1}\right)  \right\Vert _{\omega}%
^{2}+\varepsilon\left\Vert w_{0}\right\Vert ^{2}+\left\langle y_{d}%
,w_{0}\right\rangle =0\text{ .}%
\]
By Cauchy-Schwarz inequality,
\[%
\begin{array}
[c]{ll}%
\ell\left\Vert w\left(  T_{0}+T_{2}-T_{1}\right)  \right\Vert _{\omega}%
^{2}+\varepsilon\left\Vert w_{0}\right\Vert ^{2} & \leq\left\langle
Py_{d},y_{d}\right\rangle \left\langle P^{-1}w_{0},w_{0}\right\rangle \\
& \leq\left\langle Py_{d},y_{d}\right\rangle \left(  \ell\left\Vert w\left(
T_{0}+T_{2}-T_{1}\right)  \right\Vert _{\omega}^{2}+\varepsilon\left\Vert
w_{0}\right\Vert ^{2}\right)  ^{1/2}%
\end{array}
\]
where in the last line, we used $(\mathcal{O})$. Therefore, we get%
\[
\ell\left\Vert w\left(  T_{0}+T_{2}-T_{1}\right)  \right\Vert _{\omega}%
^{2}+\varepsilon\left\Vert w_{0}\right\Vert ^{2}\leq\left\langle Py_{d}%
,y_{d}\right\rangle \text{ ,}%
\]
that is,
\[
\frac{1}{\ell}\left\Vert f\right\Vert _{\omega}^{2}+\frac{1}{\varepsilon
}\left\Vert y\left(  T_{2}\right)  -y_{d}\right\Vert ^{2}\leq\left\langle
Py_{d},y_{d}\right\rangle \text{ }%
\]
where%
\[
\left\{
\begin{array}
[c]{ll}%
y^{\prime}\left(  t\right)  +Py\left(  t\right)  =0\text{ ,} & t\in\left(
T_{0},T_{2}\right)  \backslash\left\{  T_{1}\right\}  \ \text{,}\\
y\left(  T_{0}\right)  =0\text{ ,} & \\
y\left(  T_{1}\right)  =y\left(  T_{1-}\right)  +1_{\omega}\left(  -\ell
w\left(  T_{0}+T_{2}-t\right)  \right)  \text{ ,} & \\
w^{\prime}\left(  t\right)  +Pw\left(  t\right)  =0\text{ ,} & t\in\left(
T_{0},T_{2}\right)  \ \text{,}\\
w\left(  T_{0}\right)  =\frac{1}{\varepsilon}\left(  y\left(  T_{2}\right)
-y_{d}\right)  \text{ .} &
\end{array}
\right.
\]
This completes the proof. \hfill$\Box$

\bigskip

\bigskip

\subsection{Approximate impulse control}

\bigskip

Direct applications of Theorem 3.1 and Theorem 3.2, Theorem 3.3 are given now
(see \cite{Vo}\ for applications to inverse source problem). Recall that
$\omega$ is an open and nonempty subset of $\Omega$.

\bigskip

Corollary 3.1 .- \textit{Let }$0<L<T$ and $\varepsilon>0$\textit{. If one of
the statement of Theorem 3.1 holds then for any }$y_{e}\in L^{2}\left(
\Omega\right)  $\textit{, there is }$f\in L^{2}\left(  \omega\right)
$\textit{ such that the solution} $y$\textit{ to}
\[
\left\{
\begin{array}
[c]{ll}%
y^{\prime}\left(  t\right)  +Py\left(  t\right)  =0\text{ ,} & t\in\left(
0,T\right)  \backslash\left\{  L\right\}  \ \text{,}\\
y\left(  0\right)  =y_{e}\text{ ,} & \\
y\left(  L\right)  =y\left(  L_{-}\right)  +1_{\omega}f\text{ ,} &
\end{array}
\right.
\]
\textit{satisfies }%
\[
\left\Vert y\left(  T\right)  \right\Vert ^{2}\leq\varepsilon\left\Vert
y_{e}\right\Vert ^{2}\text{ \textit{and} }\left\Vert f\right\Vert _{\omega
}^{2}\leq e^{C_{3}\left(  1+\left(  \frac{1}{T-L}\right)  ^{\frac{\sigma
}{1-\sigma}}\right)  }e^{\left(  \frac{C_{3}}{T-L}\text{ln}\left(  e+\frac
{1}{\varepsilon}\right)  \right)  ^{\sigma}}\left\Vert y_{e}\right\Vert
^{2}\text{ .}%
\]

\bigskip

Proof .- We apply Theorem 3.2 with $\ell=p_{\sigma}\left(  T-L,\varepsilon
\right)  =e^{C_{3}\left(  1+\left(  \frac{1}{T-L}\right)  ^{\frac{\sigma
}{1-\sigma}}\right)  }e^{\left(  \frac{C_{3}}{T-L}\text{ln}\left(  e+\frac
{1}{\varepsilon}\right)  \right)  ^{\sigma}}$ given by Theorem 3.1 and
$T_{0}=0$, $T_{1}=L$, $T_{2}=T$ (knowing that $\left\Vert u\left(  T\right)
\right\Vert \leq\left\Vert u\left(  T-L\right)  \right\Vert $).\hfill$\Box$

\bigskip

Corollary 3.2 .- \textit{Let }$0<L<T$ and $\varepsilon>0$\textit{. If one of
the statement of Theorem 3.1 holds then for any }$y_{d}\in L^{2}\left(
\Omega\right)  $\textit{ such that }$\left\langle Py_{d},y_{d}\right\rangle
<+\infty$\textit{, there is }$f\in L^{2}\left(  \omega\right)  $\textit{ such
that the solution} $y$\textit{ to}
\[
\left\{
\begin{array}
[c]{ll}%
y^{\prime}\left(  t\right)  +Py\left(  t\right)  =0\text{ ,} & t\in\left(
0,T\right)  \backslash\left\{  L\right\}  \ \text{,}\\
y\left(  0\right)  =0\text{ ,} & \\
y\left(  L\right)  =y\left(  L_{-}\right)  +1_{\omega}f\text{ ,} &
\end{array}
\right.
\]
\textit{satisfies }$\left\Vert y\left(  T\right)  -y_{d}\right\Vert ^{2}%
\leq\varepsilon\left\langle Py_{d},y_{d}\right\rangle $\textit{ and} \textit{
}%
\[
\left\Vert f\right\Vert _{\omega}^{2}\leq\frac{1}{\lambda_{1}}e^{2C_{4}\left(
1+\left(  \frac{1}{T-L}\right)  ^{\frac{\sigma}{1-\sigma}}\right)
}e^{2\left[  \frac{T}{\varepsilon}+\left(  \frac{C_{4}}{T-L}\text{ln}\left(
\sqrt{\frac{1}{\varepsilon\lambda_{1}}}e^{\frac{2T}{\varepsilon}}\right)
\right)  ^{\sigma}\right]  }\left\langle Py_{d},y_{d}\right\rangle \text{ .}%
\]

Proof .- Recall that $\frac{d}{dt}\left\langle P^{-1}u\left(  t\right)
,u\left(  t\right)  \right\rangle +2\left\Vert u\left(  t\right)  \right\Vert
^{2}=0$ and it can be written as
\[
\frac{1}{2}\frac{d}{dt}\left\langle P^{-1}u,u\right\rangle +N\left(  t\right)
\left\langle P^{-1}u,u\right\rangle =0\text{ with }N\left(  t\right)
=\frac{\left\Vert u\left(  t\right)  \right\Vert ^{2}}{\left\langle
P^{-1}u\left(  t\right)  ,u\left(  t\right)  \right\rangle }\text{ .}%
\]
In the spirit of \cite{BT} (see also \cite{Ph}), one can check that
$N^{\prime}\left(  t\right)  \leq0$ by using Cauchy-Schwarz inequality:
$\left\Vert u\right\Vert ^{2}\leq\left\langle P^{-1}u,u\right\rangle
\left\langle Pu,u\right\rangle $ and $\frac{d}{dt}\left\Vert u\right\Vert
^{2}+2\left\langle Pu,u\right\rangle =0$. Therefore,
\[
\left\langle P^{-1}u\left(  0\right)  ,u\left(  0\right)  \right\rangle \leq
e^{2N\left(  0\right)  T}\left\langle P^{-1}u\left(  T\right)  ,u\left(
T\right)  \right\rangle \leq\frac{1}{\lambda_{1}}e^{2N\left(  0\right)
T}\left\Vert u\left(  T\right)  \right\Vert ^{2}\text{ .}%
\]
But by Theorem 3.1, it holds
\[
\left\Vert u\left(  T\right)  \right\Vert \leq e^{C_{4}\left(  1+\left(
\frac{1}{T}\right)  ^{\frac{\sigma}{1-\sigma}}\right)  }e^{\left(  \frac
{C_{4}}{T}\text{ln}\left(  \frac{\left\Vert u\left(  0\right)  \right\Vert
}{\left\Vert u\left(  T\right)  \right\Vert }\right)  \right)  ^{\sigma}%
}\left\Vert u\left(  T\right)  \right\Vert _{\omega}\text{ .}%
\]
Therefore,
\[
\left\langle P^{-1}u\left(  0\right)  ,u\left(  0\right)  \right\rangle
\leq\frac{1}{\lambda_{1}}e^{2N\left(  0\right)  T}e^{2C_{4}\left(  1+\left(
\frac{1}{T}\right)  ^{\frac{\sigma}{1-\sigma}}\right)  }e^{2\left(
\frac{C_{4}}{T}\text{ln}\left(  \frac{\left\Vert u\left(  0\right)
\right\Vert }{\left\Vert u\left(  T\right)  \right\Vert }\right)  \right)
^{\sigma}}\left\Vert u\left(  T\right)  \right\Vert _{\omega}^{2}%
\]
which implies, using $\frac{\left\Vert u\left(  0\right)  \right\Vert
}{\left\Vert u\left(  T\right)  \right\Vert }\leq\frac{\left\Vert u\left(
0\right)  \right\Vert }{\sqrt{\lambda_{1}}\left\langle P^{-1}u\left(
T\right)  ,u\left(  T\right)  \right\rangle }\leq\sqrt{\frac{N\left(
0\right)  }{\lambda_{1}}}e^{2N\left(  0\right)  T}$, the following estimate
\[
\left\langle P^{-1}u\left(  0\right)  ,u\left(  0\right)  \right\rangle
\leq\frac{1}{\lambda_{1}}e^{2N\left(  0\right)  T}e^{2C_{4}\left(  1+\left(
\frac{1}{T}\right)  ^{\frac{\sigma}{1-\sigma}}\right)  }e^{2\left(
\frac{C_{4}}{T}\text{ln}\left(  \sqrt{\frac{N\left(  0\right)  }{\lambda_{1}}%
}e^{2N\left(  0\right)  T}\right)  \right)  ^{\sigma}}\left\Vert u\left(
T\right)  \right\Vert _{\omega}^{2}\text{ . }%
\]
One conclude by distinguishing the case $N\left(  0\right)  \leq1/\varepsilon$
and the case $N\left(  0\right)  >1/\varepsilon$, that for any $\varepsilon
,T>0$,
\[
\left\langle P^{-1}u\left(  0\right)  ,u\left(  0\right)  \right\rangle
\leq\frac{1}{\lambda_{1}}e^{\frac{2T}{\varepsilon}}e^{2C_{4}\left(  1+\left(
\frac{1}{T}\right)  ^{\frac{\sigma}{1-\sigma}}\right)  }e^{2\left(
\frac{C_{4}}{T}\text{ln}\left(  \sqrt{\frac{1}{\varepsilon\lambda_{1}}%
}e^{\frac{2T}{\varepsilon}}\right)  \right)  ^{\sigma}}\left\Vert u\left(
T\right)  \right\Vert _{\omega}^{2}+\varepsilon\left\Vert u\left(  0\right)
\right\Vert ^{2}\text{. }%
\]
It remains to apply Theorem 3.3 with $\ell=\frac{1}{\lambda_{1}}e^{\frac
{2T}{\varepsilon}}e^{2C_{4}\left(  1+\left(  \frac{1}{T-L}\right)
^{\frac{\sigma}{1-\sigma}}\right)  }e^{2\left(  \frac{C_{4}}{T-L}%
\text{ln}\left(  \sqrt{\frac{1}{\varepsilon\lambda_{1}}}e^{\frac
{2T}{\varepsilon}}\right)  \right)  ^{\sigma}}$ and $T_{0}=0$, $T_{1}=L$,
$T_{2}=T$.\hfill$\Box$

\bigskip

\bigskip

\subsection{Null controllability with measurable set in time}

\bigskip

Recall that $\omega$ is an open and nonempty subset of $\Omega$.

\bigskip

Theorem 3.4 .- \textit{Let }$T>0$\textit{ and }$E\subset\left(  0,T\right)
$\textit{ a set of positive measure. If one of the statement of Theorem 3.1
holds then for any }$y^{0}\in L^{2}\left(  \Omega\right)  $\textit{, there is
}$f\in L^{2}\left(  \omega\times E\right)  $\textit{ such that the solution
}$y$\textit{ to}
\[
\left\{
\begin{array}
[c]{ll}%
y^{\prime}\left(  t\right)  +Py\left(  t\right)  =1_{\omega\times E}f\text{ ,}
& t\in\left(  0,T\right)  \ \text{,}\\
y\left(  0\right)  =y^{0}\text{ ,} &
\end{array}
\right.
\]
\textit{satisfies }$y\left(  T\right)  =0$.

\bigskip

Proof .- The proof is divided into three steps.

\bigskip

\textit{Step 1: Observability estimate with measurable set in time. }Based on
a telescoping series method (see \cite{Mi}, \cite{Mi2} and already exploited
in \cite{PW}, \cite{PWZ}, \cite{AEWZ}, \cite{EMZ}, \cite{Z}, \cite{WZ},
\cite{LiZ}, \cite{YZ}, \cite{Ph2}), the statement $\left(  ii\right)  $ in
Theorem 3.1 implies the following observability: The solution $u$ to
\[
\left\{
\begin{array}
[c]{ll}%
u^{\prime}\left(  t\right)  +Pu\left(  t\right)  =0\text{ ,} & t\in\left(
0,T\right)  \ \text{,}\\
u\left(  0\right)  \in L^{2}\left(  \Omega\right)  \text{ ,} &
\end{array}
\right.
\]
satisfies%
\[
\left\Vert u\left(  T\right)  \right\Vert ^{2}\leq K\int_{E}\left\Vert
u\left(  T-t\right)  \right\Vert _{\omega}^{2}dt\text{ .}%
\]
Here, $K$ is a constant only depending on $\left(  P,\Omega,\omega
,\sigma,\left\vert E\right\vert \right)  $. Further, if $E=\left(  0,T\right)
$, then $K=C$exp$\left(  \frac{C}{T^{\frac{\sigma}{1-\sigma}}}\right)  $ for
some $C=C\left(  P,\Omega,\omega,\sigma\right)  $.

\bigskip

\textit{Step 2: Approximate controllability. }Let $\varepsilon>0$. Consider
the functional $J_{\varepsilon}$ defined on $L^{2}\left(  \Omega\right)  $
given by%
\[
J_{\varepsilon}\left(  u_{0}\right)  =\frac{K}{2}\int_{E}\left\Vert u\left(
T-t\right)  \right\Vert _{\omega}^{2}dt+\frac{\varepsilon}{2}\left\Vert
u_{0}\right\Vert ^{2}-\left\langle y^{0},u\left(  T\right)  \right\rangle
\text{ ,}%
\]
where%
\[
\left\{
\begin{array}
[c]{ll}%
u^{\prime}\left(  t\right)  +Pu\left(  t\right)  =0\text{ ,} & t\in\left(
0,T\right)  \ \text{,}\\
u\left(  0\right)  =u_{0}\text{ .} &
\end{array}
\right.
\]
Notice that $J_{\varepsilon}$ is strictly convex, $C^{1}$ and coercive and
therefore $J_{\varepsilon}$ has a unique minimizer $w_{\varepsilon,0}\in
L^{2}\left(  \Omega\right)  $, i.e. $J_{\varepsilon}(w_{\varepsilon
,0})=\underset{u_{0}\in L^{2}\left(  \Omega\right)  }{\text{min}%
}J_{\varepsilon}(u_{0})$. Set%
\[
\left\{
\begin{array}
[c]{ll}%
w_{\varepsilon}^{\prime}\left(  t\right)  +Pw_{\varepsilon}\left(  t\right)
=0\text{ ,} & t\in\left(  0,T\right)  \ \text{,}\\
w_{\varepsilon}\left(  0\right)  =w_{\varepsilon,0}\text{ ,} &
\end{array}
\right.  \text{ and }\left\{
\begin{array}
[c]{ll}%
h^{\prime}\left(  t\right)  +Ph\left(  t\right)  =0\text{ ,} & t\in\left(
0,T\right)  \ \text{,}\\
h\left(  0\right)  =h_{0}\text{ .} &
\end{array}
\right.
\]
Since $J_{\varepsilon}^{\prime}(w_{\varepsilon,0})h_{0}=0$ for any $h_{0}\in
L^{2}\left(  \Omega\right)  $, we have%
\[
K\int_{E}\left\langle w_{\varepsilon}\left(  T-t\right)  ,h\left(  T-t\right)
\right\rangle _{\omega}dt+\varepsilon\left\langle w_{\varepsilon,0}%
,h_{0}\right\rangle -\left\langle y^{0},h\left(  T\right)  \right\rangle
=0\quad\forall h_{0}\in L^{2}\left(  \Omega\right)  \text{ .}%
\]
But the solution $y_{\varepsilon}$ to
\[
\left\{
\begin{array}
[c]{ll}%
y_{\varepsilon}^{\prime}\left(  t\right)  +Py_{\varepsilon}\left(  t\right)
=1_{\omega\times E}f_{\varepsilon}\text{ ,} & t\in\left(  0,T\right)
\ \text{,}\\
y_{\varepsilon}\left(  0\right)  =y^{0}\text{ ,} &
\end{array}
\right.
\]
satisfies
\[
\left\langle y_{\varepsilon}\left(  T\right)  ,u\left(  0\right)
\right\rangle -\left\langle y^{0},u\left(  T\right)  \right\rangle =\int
_{E}\left\langle f_{\varepsilon}\left(  \cdot,t\right)  ,u\left(  T-t\right)
\right\rangle _{\omega}dt\quad\forall u\left(  0\right)  \in L^{2}\left(
\Omega\right)
\]
which means%
\[
-\int_{E}\left\langle f_{\varepsilon}\left(  \cdot,t\right)  ,h\left(
T-t\right)  \right\rangle _{\omega}dt+\left\langle y_{\varepsilon}\left(
T\right)  ,h_{0}\right\rangle -\left\langle y^{0},h\left(  T\right)
\right\rangle =0\quad\forall h_{0}\in L^{2}\left(  \Omega\right)  \text{ .}%
\]
By choosing $f_{\varepsilon}\left(  \cdot,t\right)  =-Kw_{\varepsilon}\left(
T-t\right)  $, we deduce that the solution $y_{\varepsilon}$ satisfies
\[
\varepsilon w_{\varepsilon,0}=y_{\varepsilon}\left(  T\right)  \text{ .}%
\]
Further,
\[
K\int_{E}\left\Vert w_{\varepsilon}\left(  T-t\right)  \right\Vert _{\omega
}^{2}dt+\varepsilon\left\Vert w_{\varepsilon,0}\right\Vert ^{2}=\frac{1}%
{K}\int_{E}\left\Vert f_{\varepsilon}\left(  \cdot,t\right)  \right\Vert
_{\omega}^{2}dt+\frac{1}{\varepsilon}\left\Vert y_{\varepsilon}\left(
T\right)  \right\Vert ^{2}\text{ .}%
\]
Moreover, taking $h_{0}=w_{\varepsilon,0}$ into $J_{\varepsilon}^{\prime
}(w_{\varepsilon,0})h_{0}=0$, we get%
\[
K\int_{E}\left\Vert w_{\varepsilon}\left(  T-t\right)  \right\Vert _{\omega
}^{2}dt+\varepsilon\left\Vert w_{\varepsilon,0}\right\Vert ^{2}-\left\langle
y^{0},w_{\varepsilon}\left(  T\right)  \right\rangle =0\text{ .}%
\]
By Cauchy-Schwarz inequality,
\[%
\begin{array}
[c]{ll}%
K\displaystyle\int_{E}\left\Vert w_{\varepsilon}\left(  T-t\right)
\right\Vert _{\omega}^{2}dt+\varepsilon\left\Vert w_{\varepsilon,0}\right\Vert
^{2} & \leq\left\Vert y^{0}\right\Vert \left\Vert w_{\varepsilon}\left(
T\right)  \right\Vert \\
& \leq\left\Vert y^{0}\right\Vert \left(  K\displaystyle\int_{E}\left\Vert
w_{\varepsilon}\left(  T-t\right)  \right\Vert _{\omega}^{2}dt\right)  ^{1/2}%
\end{array}
\]
where in the last line, we used the observability estimate with measurable set
in time. Therefore, we get%
\[
K\int_{E}\left\Vert w_{\varepsilon}\left(  T-t\right)  \right\Vert _{\omega
}^{2}dt+2\varepsilon\left\Vert w_{\varepsilon,0}\right\Vert ^{2}\leq\left\Vert
y^{0}\right\Vert ^{2}\text{ ,}%
\]
that is,
\[
\frac{1}{K}\int_{E}\left\Vert f_{\varepsilon}\left(  \cdot,t\right)
\right\Vert _{\omega}^{2}dt+\frac{2}{\varepsilon}\left\Vert y_{\varepsilon
}\left(  T\right)  \right\Vert ^{2}\leq\left\Vert y^{0}\right\Vert ^{2}\text{
}%
\]
where%
\[
\left\{
\begin{array}
[c]{ll}%
y_{\varepsilon}^{\prime}\left(  t\right)  +Py_{\varepsilon}\left(  t\right)
=1_{\omega\times E}f_{\varepsilon}\text{ ,} & t\in\left(  0,T\right)
\ \text{,}\\
y_{\varepsilon}\left(  0\right)  =y^{0}\text{ ,} & \\
f_{\varepsilon}\left(  x,t\right)  =-Kw_{\varepsilon}\left(  x,T-t\right)
\text{ ,} & \left(  x,t\right)  \in\Omega\times\left(  0,T\right)
\ \text{,}\\
w_{\varepsilon}^{\prime}\left(  t\right)  +Pw_{\varepsilon}\left(  t\right)
=0\text{ ,} & t\in\left(  0,T\right)  \ \text{,}\\
w_{\varepsilon}\left(  T\right)  =\frac{1}{\varepsilon}y_{\varepsilon}\left(
T\right)  \text{ .} &
\end{array}
\right.
\]

\bigskip

\textit{Step 3: Convergence of the control function.} We refer to
\cite[p.571]{Zu}. Since $w_{\varepsilon}\left(  T-\cdot\right)  $ is bounded
in $L^{2}\left(  \omega\times E\right)  $ and $\sqrt{\varepsilon
}w_{\varepsilon,0}$ is bounded in $L^{2}\left(  \Omega\right)  $, one can
deduce that, for some function $w\left(  T-\cdot\right)  $ in $L^{2}\left(
\omega\times E\right)  $, $w_{\varepsilon}\left(  T-\cdot\right)  $ weakly
converge to $w\left(  T-\cdot\right)  $ in $L^{2}\left(  \omega\times
E\right)  $ and $\varepsilon w_{\varepsilon,0}$ tends to zero in $L^{2}\left(
\Omega\right)  $. Therefore the identity
\[
K\int_{E}\left\langle w_{\varepsilon}\left(  T-t\right)  ,h\left(  T-t\right)
\right\rangle _{\omega}dt+\varepsilon\left\langle w_{\varepsilon,0}%
,h_{0}\right\rangle -\left\langle y^{0},h\left(  T\right)  \right\rangle
=0\quad\forall h_{0}\in L^{2}\left(  \Omega\right)  \text{ ,}%
\]
becomes when $\varepsilon\rightarrow0$, as%
\[
K\int_{E}\left\langle w\left(  T-t\right)  ,h\left(  T-t\right)  \right\rangle
_{\omega}dt-\left\langle y^{0},h\left(  T\right)  \right\rangle =0\quad\forall
h_{0}\in L^{2}\left(  \Omega\right)  \text{ .}%
\]
But the solution $y$ to
\[
\left\{
\begin{array}
[c]{ll}%
y^{\prime}\left(  t\right)  +Py\left(  t\right)  =1_{\omega\times E}f\text{ ,}
& t\in\left(  0,T\right)  \ \text{,}\\
y\left(  0\right)  =y^{0}\text{ ,} &
\end{array}
\right.
\]
satisfies
\[
-\int_{E}\left\langle f\left(  \cdot,t\right)  ,h\left(  T-t\right)
\right\rangle _{\omega}dt+\left\langle y\left(  T\right)  ,h_{0}\right\rangle
-\left\langle y^{0},h\left(  T\right)  \right\rangle =0\quad\forall h_{0}\in
L^{2}\left(  \Omega\right)  \text{ .}%
\]
By choosing $f\left(  \cdot,t\right)  =-Kw\left(  T-t\right)  $, it follows
that the solution $y$ satisfies
\[
y\left(  T\right)  =0\text{ .}%
\]
This completes the proof. \hfill$\Box$

\bigskip

\bigskip

\subsection{Finite time stabilization}

\bigskip

Recall that $\left\Vert \cdot\right\Vert $ and $\left\langle \cdot
,\cdot\right\rangle $ are the norm and the inner product of $L^{2}\left(
\Omega\right)  $ respectively.

\bigskip

Assume that there are two positive constants $c=c\left(  \Omega\right)  $ and
$\rho=\rho\left(  d\right)  $ such that
\[
\text{Card}\left\{  \lambda_{i}\leq\Lambda\right\}  =\sum\limits_{\lambda
_{i}\leq\Lambda}1\leq c\Lambda^{1/\rho}\text{ .}%
\]
Such estimate can be provided by the Weyl asymptotic formula $\lambda_{k}\sim
C\left(  \Omega\right)  k^{\rho}$ as $k\rightarrow\infty$. In particular, if
$P=-\Delta$, then $\rho=\frac{2}{d}$; And if $P=\Delta^{2}$, then $\rho
=\frac{4}{d}$. In the case of the one-dimensional degenerate operator
$P=\mathcal{P}$, we have $\rho=2$.

\bigskip

Define an increasing sequence $\left(  t_{m}\right)  _{m\geq0}$ converging to
$T>0$ by
\[
t_{m}=T\left(  1-\frac{1}{b^{m}}\right)  \text{ for some }b>1\text{ .}%
\]
Introduce a linear bounded operator $\mathcal{F}_{m}$ from $L^{2}\left(
\Omega\right)  $ into $L^{2}\left(  \omega\right)  $ in the following manner:%
\[%
\begin{array}
[c]{cccl}%
\mathcal{F}_{m}: & L^{2}\left(  \Omega\right)  & \rightarrow & L^{2}\left(
\omega\right) \\
& \vartheta & \mapsto & \sum\limits_{\lambda_{j}\leq\Lambda_{m}}\left\langle
\vartheta,\Phi_{j}\right\rangle f_{j}%
\end{array}
\]
where
\[
\Lambda_{m}:=\lambda_{1}+\left(  \frac{\eta}{T}\frac{b}{b-1}\right)
b^{\left(  \beta+1\right)  m}\text{ with }\eta>1\text{ ,}\quad\beta
:=\frac{\sigma}{1-\sigma}\text{ ,}%
\]
and $f_{j}$ is the impulse control of the heat equation associated to the
eigenfunction $\Phi_{j}$ (see Corollary 3.1):
\[
\left\{
\begin{array}
[c]{ll}%
y_{j}^{\prime}\left(  t\right)  +Py_{j}\left(  t\right)  =0\text{ ,} &
t\in\left(  t_{m},t_{m+1}\right)  \backslash\left\{  \frac{t_{m}+t_{m+1}}%
{2}\right\}  \ \text{,}\\
y_{j}\left(  t_{m}\right)  =\Phi_{j}\text{ ,} & \\
y_{j}\left(  \frac{t_{m}+t_{m+1}}{2}\right)  =y_{j}\left(  \left(  \frac
{t_{m}+t_{m+1}}{2}\right)  _{-}\right)  +1_{\omega}f_{j}\text{ ,} &
\end{array}
\right.
\]
satisfying\textit{ }%
\[
\left\Vert y_{j}\left(  t_{m+1}\right)  \right\Vert ^{2}\leq\frac{e^{-\eta
b^{\beta m}}}{\sum\limits_{\lambda_{i}\leq\Lambda_{m}}1}\text{ and }\left\Vert
f_{j}\right\Vert _{\omega}^{2}\leq e^{C_{3}\left(  1+\left(  \frac{2}%
{t_{m+1}-t_{m}}\right)  ^{\beta}\right)  }e^{\left(  \frac{2C_{3}}%
{t_{m+1}-t_{m}}\text{ln}\left(  e+e^{\eta b^{\beta m}}\sum\limits_{\lambda
_{i}\leq\Lambda_{m}}1\right)  \right)  ^{\sigma}}\text{ .}%
\]
Here, $C_{3}>0$ and $\sigma\in\left(  0,1\right)  $ are the constants given in
Theorem 3.1. Notice that
\[
\left\Vert \mathcal{F}_{m}\right\Vert _{L^{2}\left(  \Omega\right)
\rightarrow L^{2}\left(  \omega\right)  }^{2}\leq\sum\limits_{\lambda_{j}%
\leq\Lambda_{m}}\left\Vert f_{j}\right\Vert _{\omega}^{2}\text{ .}%
\]

\bigskip

\bigskip

Theorem 3.5 .- \textit{Let }$\omega$\textit{ be an open and nonempty subset of
}$\Omega$\textit{. Suppose that one of the statement of Theorem 3.1 holds and}%
\[
\text{Card}\left\{  \lambda_{i}\leq\Lambda\right\}  \leq c\Lambda^{1/\rho
}\text{ }\quad\text{\ for any }\Lambda>0\text{ .}%
\]
\textit{Then, for any }$T>0$ \textit{there are }$b,\eta>1$\textit{ and
}$C,K>0$\textit{ such that for any }$z_{0}\in L^{2}\left(  \Omega\right)
$\textit{, the solution }$z$\textit{ to }%
\[
\left\{
\begin{array}
[c]{ll}%
z^{\prime}\left(  t\right)  +Pz\left(  t\right)  =0\text{ ,} & t\in
\mathbb{R}^{+}\backslash\bigcup\limits_{m\geq0}\left(  \frac{t_{m}+t_{m+1}}%
{2}\right)  \ \text{,}\\
z\left(  \frac{t_{m}+t_{m+1}}{2}\right)  =z\left(  \left(  \frac{t_{m}%
+t_{m+1}}{2}\right)  _{-}\right)  +1_{\omega}\mathcal{F}_{m}\left(  z\left(
t_{m}\right)  \right)  \text{ ,} & \text{\textit{for any integer} }%
m\geq0\text{ ,}\\
z\left(  0\right)  =z_{0}\text{ ,} &
\end{array}
\right.
\]
\textit{satisfies }$\left\Vert z\left(  t\right)  \right\Vert \leq
Ce^{-\frac{1}{K}\left(  \frac{T}{T-t}\right)  ^{\frac{\sigma}{1-\sigma}}%
}\left\Vert z_{0}\right\Vert $\textit{ for any }$0\leq t<T_{-}$\textit{.
Further, }$\underset{m\rightarrow\infty}{\text{lim}}\left\Vert \mathcal{F}%
_{m}\left(  z\left(  t_{m}\right)  \right)  \right\Vert =0$.

\bigskip

Proof .- We start to focus on the solution $z$ on interval $\left(
t_{m},t_{m+1}\right)  $ with initial data $z\left(  t_{m}\right)
=\sum\limits_{j\geq1}a_{j}\Phi_{j}$ in $L^{2}\left(  \Omega\right)  $.
Introduce the initial datum $\phi\left(  t_{m}\right)  =\sum\limits_{\lambda
_{j}>\Lambda_{m}}a_{j}\Phi_{j}$ and $\psi\left(  t_{m}\right)  =\sum
\limits_{\lambda_{j}\leq\Lambda_{m}}a_{j}\Phi_{j}$ associated to the solution
of $\phi^{\prime}\left(  t\right)  +P\phi\left(  t\right)  =0$ and
\[
\left\{
\begin{array}
[c]{ll}%
\psi^{\prime}\left(  t\right)  +P\psi\left(  t\right)  =0\text{ ,} &
t\in\left(  t_{m},t_{m+1}\right)  \backslash\left\{  \frac{t_{m}+t_{m+1}}%
{2}\right\}  \ \text{,}\\
\psi\left(  \frac{t_{m}+t_{m+1}}{2}\right)  =\psi\left(  \left(  \frac
{t_{m}+t_{m+1}}{2}\right)  _{-}\right)  +1_{\omega}\displaystyle\sum
\limits_{\lambda_{j}\leq\Lambda_{m}}a_{j}f_{j}\text{ .} &
\end{array}
\right.
\]
Therefore, $\phi\left(  t_{m+1}\right)  =\sum\limits_{\lambda_{j}>\Lambda_{m}%
}a_{j}e^{-\lambda_{j}\left(  t_{m+1}-t_{m}\right)  }\Phi_{j}$ and%
\[
\left\Vert \phi\left(  t_{m+1}\right)  \right\Vert \leq e^{-\Lambda_{m}\left(
t_{m+1}-t_{m}\right)  }\left\Vert z\left(  t_{m}\right)  \right\Vert \text{ .}%
\]
But we have chosen $\Lambda_{m}>\lambda_{1}$ in order that $\eta b^{\beta
m}\leq\Lambda_{m}\left(  t_{m+1}-t_{m}\right)  $. It implies that
\[
\left\Vert \phi\left(  t_{m+1}\right)  \right\Vert \leq e^{-\eta b^{\beta m}%
}\left\Vert z\left(  t_{m}\right)  \right\Vert \text{ .}%
\]
On the other hand, the solution $\psi$ satisfies $\psi=\sum\limits_{\lambda
_{j}\leq\Lambda_{m}}a_{j}y_{j}$ and
\[
\left\Vert \psi\left(  t_{m+1}\right)  \right\Vert \leq\sum\limits_{\lambda
_{j}\leq\Lambda_{m}}\left\vert a_{j}\right\vert \sqrt{\frac{e^{-\eta b^{\beta
m}}}{\sum\limits_{\lambda_{i}\leq\Lambda_{m}}1}}\leq e^{-\frac{1}{2}\eta
b^{\beta m}}\left\Vert v\left(  t_{m}\right)  \right\Vert \text{ .}%
\]
Consequently, we have
\[
\left\Vert z\left(  t_{m+1}\right)  \right\Vert \leq\left\Vert \phi\left(
t_{m+1}\right)  \right\Vert +\left\Vert \psi\left(  t_{m+1}\right)
\right\Vert \leq e^{1-\frac{1}{2}\eta b^{\beta m}}\left\Vert z\left(
t_{m}\right)  \right\Vert \text{ ,}%
\]
which implies by induction that for any $m\geq1$,
\[
\left\Vert z\left(  t_{m}\right)  \right\Vert ^{2}\leq e^{2m-\eta b^{\beta m}%
}\left\Vert z\left(  t_{0}\right)  \right\Vert ^{2}\text{ .}%
\]
Now, we treat the boundeness of the control associated to $\psi$: Notice that
$\sum\limits_{\lambda_{j}\leq\Lambda_{m}}a_{j}f_{j}:=\mathcal{F}_{m}\left(
z\left(  t_{m}\right)  \right)  $, and then, by Cauchy-Schwarz and Young
inequalities,%
\[%
\begin{array}
[c]{ll}%
\left\Vert \mathcal{F}_{m}\left(  z\left(  t_{m}\right)  \right)  \right\Vert
_{\omega}^{2} & \leq\displaystyle\int_{\omega}\left(  \sum\limits_{\lambda
_{j}\leq\Lambda_{m}}\left\vert a_{j}\right\vert \left\vert f_{j}\right\vert
\right)  ^{2}\leq\displaystyle\sum\limits_{\lambda_{j}\leq\Lambda_{m}%
}\left\vert a_{j}\right\vert ^{2}\displaystyle\sum\limits_{j\leq\Lambda_{m}%
}\left\Vert f_{j}\right\Vert _{\omega}^{2}\\
& \leq\left\Vert z\left(  t_{m}\right)  \right\Vert ^{2}\displaystyle\sum
\limits_{\lambda_{j}\leq\Lambda_{m}}e^{C_{3}\left(  1+\left(  \frac{2}%
{t_{m+1}-t_{m}}\right)  ^{\beta}\right)  }e^{\left(  \frac{2C_{3}}%
{t_{m+1}-t_{m}}\text{ln}\left(  e+\sum\limits_{\lambda_{i}\leq\Lambda_{m}%
}1e^{\eta b^{\beta m}}\right)  \right)  ^{\sigma}}\\
& \leq\left\Vert z\left(  t_{m}\right)  \right\Vert ^{2}e^{C_{3}\left(
1+\left(  \frac{2}{t_{m+1}-t_{m}}\right)  ^{\beta}\right)  }\displaystyle\sum
\limits_{\lambda_{j}\leq\Lambda_{m}}e^{\left(  \frac{4C_{3}}{t_{m+1}-t_{m}%
}\right)  ^{\beta}+\frac{1}{2}\text{ln}\left(  \sum\limits_{\lambda_{i}%
\leq\Lambda_{m}}1e^{\eta b^{\beta m}}\right)  }\\
& \leq\left\Vert z\left(  t_{m}\right)  \right\Vert ^{2}e^{C_{3}\left(
1+\left(  \frac{2}{t_{m+1}-t_{m}}\right)  ^{\beta}\right)  }e^{\left(
\frac{4C_{3}}{t_{m+1}-t_{m}}\right)  ^{\frac{\sigma}{1-\sigma}}}e^{\frac{1}%
{2}\eta b^{\beta m}}\displaystyle\left(  \sum\limits_{\lambda_{i}\leq
\Lambda_{m}}1\right)  ^{\theta+1}\\
& \leq e^{2m-\eta b^{\beta m}}\left\Vert z\left(  t_{0}\right)  \right\Vert
^{2}e^{C_{3}\left(  1+\left(  \frac{2}{T}\frac{b^{m+1}}{b-1}\right)  ^{\beta
}\right)  }e^{\left(  \frac{4C_{3}}{T}\frac{b^{m+1}}{b-1}\right)  ^{\beta}%
}e^{\frac{1}{2}\eta b^{\beta m}}\left(  c\Lambda_{m}^{1/\rho}\right)
^{\theta+1}\\
& \leq e^{2m-\frac{1}{2}\eta b^{\beta m}}\left\Vert z\left(  t_{0}\right)
\right\Vert ^{2}e^{C_{3}}e^{\left(  C_{3}+\left(  2C_{3}\right)  ^{\beta
}\right)  \left(  \frac{2}{T}\frac{b}{b-1}\right)  ^{\beta}b^{\beta m}}\left(
c\left(  \lambda_{1}+\left(  \frac{\eta}{T}\frac{b}{b-1}\right)  b^{\left(
\beta+1\right)  m}\right)  \right)  ^{\frac{\theta+1}{\rho}}\text{ ,}%
\end{array}
\]
where in the last line we used the definition of $\Lambda_{m}$. Next, we
choose $\eta>1$, precisely
\[
\eta=1+4\left(  C_{3}+\left(  2C_{3}\right)  ^{\beta}\right)  \left(  \frac
{2}{T}\frac{b}{b-1}\right)  ^{\beta}%
\]
in order that $-\frac{1}{2}\eta b^{\beta m}+\left(  C_{3}+\left(  \frac{C_{3}%
}{\theta}\right)  ^{\beta}\right)  \left(  \frac{2}{T}\frac{b}{b-1}\right)
^{\beta}b^{\beta m}\leq-\frac{1}{4}\eta b^{\beta m}$.

Since $b>1$ and $b^{\left(  \beta+1\right)  \left(  \theta+1\right)  m/\rho
}\leq\left(  \frac{8\left(  \beta+1\right)  \left(  \theta+1\right)  }%
{\beta\rho\eta}\right)  ^{\frac{\left(  \beta+1\right)  \left(  \theta
+1\right)  }{\beta\rho}}e^{\frac{1}{8}\eta b^{\beta m}}$, we obtain for some
constant $C_{5}:=\left(  c\left(  \lambda_{1}+\left(  \frac{\eta}{T}\frac
{b}{b-1}\right)  \right)  \right)  ^{\frac{\theta+1}{\rho}}>0$ that for any
$m\geq1$,
\[
\left\Vert \mathcal{F}_{m}\left(  v\left(  t_{m}\right)  \right)  \right\Vert
_{\omega}^{2}\leq C_{5}e^{2m-\frac{1}{8}\eta b^{\beta m}}\left\Vert v\left(
0\right)  \right\Vert ^{2}\text{ .}%
\]
Finally, let $t\geq0$, then there is $m\geq0$ such that $t\in\left[
t_{m},t_{m+1}\right]  $. We distinguish four cases: If $t\in\left[
0,t_{1}/2\right)  $, then
\[
\left\Vert z\left(  t\right)  \right\Vert ^{2}\leq\left\Vert z\left(
0\right)  \right\Vert ^{2}\text{ ;}%
\]
If $t\in\left[  t_{1}/2,t_{1}\right)  $, then
\[
\left\Vert z\left(  t\right)  \right\Vert ^{2}\leq\left\Vert z\left(  \left(
t_{1}/2\right)  _{-}\right)  +1_{\omega}\mathcal{F}_{0}\left(  z\left(
t_{0}\right)  \right)  \right\Vert ^{2}\leq2\left(  1+\left\Vert
\mathcal{F}_{0}\right\Vert ^{2}\right)  \left\Vert z\left(  0\right)
\right\Vert ^{2}\text{ .}%
\]
If $t\in\left[  t_{m},\frac{t_{m}+t_{m+1}}{2}\right)  $ and $m\geq1$, then
\[
\left\Vert z\left(  t\right)  \right\Vert ^{2}\leq\left\Vert z\left(
t_{m}\right)  \right\Vert ^{2}\leq e^{2m-\eta b^{\beta m}}\left\Vert z\left(
0\right)  \right\Vert ^{2}\text{ ;}%
\]
If $t\in\left[  \frac{t_{m}+t_{m+1}}{2},t_{m+1}\right)  $ and $m\geq1$, then
\[%
\begin{array}
[c]{ll}%
\left\Vert z\left(  t\right)  \right\Vert ^{2} & \leq\left\Vert v\left(
\left(  \frac{t_{m}+t_{m+1}}{2}\right)  _{-}\right)  +1_{\omega}%
\mathcal{F}_{m}\left(  z\left(  t_{m}\right)  \right)  \right\Vert ^{2}\\
& \leq2\left(  1+C_{5}\right)  b^{\left(  \beta+1\right)  \left(
\theta+1\right)  m/\rho}e^{2m-\frac{1}{8}\eta b^{\beta m}}\left\Vert z\left(
0\right)  \right\Vert ^{2}\text{ .}%
\end{array}
\]
Consequently, for any $t\in\left[  t_{m},t_{m+1}\right]  $, it holds%
\[
b^{\beta m}\leq\left(  \frac{T}{T-t}\right)  ^{\beta}\leq b^{\beta m+\beta}%
\]
and
\[
\left\Vert z\left(  t\right)  \right\Vert ^{2}\leq2\left(  1+C_{5}+\left\Vert
\mathcal{F}_{0}\right\Vert ^{2}\right)  e^{-\frac{1}{16}\eta b^{\beta m}%
}\left\Vert z\left(  0\right)  \right\Vert ^{2}%
\]
by choosing $b=e^{32/\left(  \beta\eta\right)  }$. One conclude that
$e^{-\frac{1}{16}\eta b^{\beta m}}\leq e^{-\frac{1}{16}\eta\left(  \frac{T}%
{b}\frac{1}{T-t}\right)  ^{\beta}}$ and
\[
\left\Vert z\left(  t\right)  \right\Vert ^{2}\leq2\left(  1+C_{5}+\left\Vert
\mathcal{F}_{0}\right\Vert ^{2}\right)  e^{-\frac{1}{16}\eta e^{32/\eta
}\left(  \frac{T}{T-t}\right)  ^{\beta}}\left\Vert z\left(  0\right)
\right\Vert ^{2}\text{ .}%
\]
This completes the proof. \hfill$\Box$

\bigskip

\bigskip

\bigskip
\end{document}